\documentclass[11pt,a4paper]{article}

\usepackage[a4paper,margin=30mm]{geometry}
\usepackage[T1]{fontenc}
\usepackage[utf8]{inputenc}
\usepackage{mathpazo}
\usepackage{courier}

\usepackage{microtype}
\usepackage{amsmath,amssymb,amsthm,amscd,mathtools}
\usepackage{mathrsfs}
\usepackage{esint}
\usepackage{enumitem}
\usepackage{xcolor}
\usepackage{hyperref}

\usepackage{tikz-cd}
\hypersetup{
	colorlinks=true,
	linkcolor=blue!55!black,
	citecolor=blue!55!black,
	urlcolor=blue!55!black,
	pdftitle={Liouville systems I: Blow-up analysis and orthogonal decomposition},
	pdfauthor={Mohameden Ahmedou, Jonas Fey, Zhengni Hu, Lei Zhang},
	pdfsubject={Blow-up solutions, amplitude variables, and orthogonal decomposition}
}
\usepackage{cite,needspace}
\allowdisplaybreaks

\numberwithin{equation}{section}

\theoremstyle{plain}
\newtheorem{theorem}{Theorem}[section]
\newtheorem{lemma}{Lemma}[section]
\newtheorem{proposition}{Proposition}[section]
\newtheorem{corollary}{Corollary}[section]

\theoremstyle{definition}
\newtheorem{definition}{Definition}[section]

\theoremstyle{remark}
\newtheorem{remark}{Remark}[section]

\newcommand{\R}{\mathbb R}

\newcommand{\cH}{\mathcal H}
\newcommand{\cO}{\mathcal O}
\newcommand{\la}{\langle}
\newcommand{\ra}{\rangle}
\newcommand{\oH}{\mathring{\mathrm H}^{1}}
\newcommand{\ff}{\mathcal F^N_{\rho^*,h}}
\newcommand{\eps}{\varepsilon}
\newcommand{\lam}{\lambda}
\newcommand{\rg}{\mathfrak g}
\newcommand{\rf}{\mathrm f}
\newcommand{\fJ}{\mathfrak c}
\title{\bfseries Liouville systems I:\\Blow-up analysis and orthogonal decomposition}
\author{Mohameden Ahmedou\thanks{Mathematisches Institut, Justus-Liebig-Universit\"at Giessen, Germany.}
	\and Jonas Fey\thanks{Corresponding author.
        Mathematisches Institut, Justus-Liebig-Universit\"at Giessen, Germany.
        Email: \href{mailto:jonas.fey@math.uni-giessen.de}{\texttt{jonas.fey@math.uni-giessen.de}}.}
	\and Zhengni Hu\thanks{Department of Mathematics, Shanghai Jiao Tong University, Shanghai, China.}
	\and Lei Zhang\thanks{Department of Mathematics, University of Florida, Gainesville, FL, USA.}}
\date{November 4, 2026}
\begin{document}
	\maketitle
	\begin{abstract}
		We study blow-up solutions of generalized Liouville systems on a closed
		surface at and near critical parameter hypersurfaces. Under a Morse
		assumption on the reduced functional, we establish the following results.\par 
		First, we introduce amplitude variables in the approximate manifold and
		prove that every sequence of blow-up solutions at a fixed critical
		parameter, or along an admissible sequence from the lower side, admits
		a unique minimizing decomposition up to permutation. The remainder
		converges to zero and is orthogonal to the
		center, scale, and amplitude directions. The concentration scales are
		comparable, the limiting configuration satisfies the scale balance
		conditions, and the amplitudes vanish after multiplication by the
		logarithmic scale. Second, we determine the leading parameter defect
		in both critical decay regimes and prove that, along each fixed
		admissible path, all sufficiently large solutions lie in finitely many
		disjoint bubbling neighborhoods.  
		Finally, with an extra assumption that the leading coefficient
		is nonzero at every balanced critical configuration, the solution set
		at the fixed critical parameter is compact.
		
		The results provide the analytic foundation for
		the study of resonant degree by critical points at infinity.
	\end{abstract}
	\medskip
	\noindent\textbf{Keywords.} Liouville systems; critical parameters; blow-up analysis; orthogonal decomposition.
	
	\noindent\textbf{2020 Mathematics Subject Classification.}
	35J60, 35B44, 35J50.
	
	\section{Introduction}
	
	Let $(M,g)$ be a closed Riemann surface with $\operatorname{Vol}_g(M)=1$.
	We consider the generalized Liouville system
	\begin{equation}\label{I:eq:L-System}
		-\Delta_g u_i=\sum_{j=1}^n a_{ij}\rho_j
		\left(\frac{h_je^{u_j}}{\int_M h_je^{u_j}\,dV_g}-1\right),
		\qquad i\in I:=\{1,\ldots,n\},
	\end{equation}
	where $h_i\in C^3(M)$ is positive and $\rho_i>0$. 
	Denote  $A=(a_{ij}), A^{-1}=(a^{ij})$. 
	Throughout the paper,
	we assume
	\begin{enumerate}[label=(H\arabic*),ref=(H\arabic*)]
		\item\label{I:assump:H1} $A=(a_{ij})$ is symmetric, nonnegative,
		irreducible, and invertible;
		\item\label{I:assump:H2} 
		$
		a^{ii}\le0, a^{ij}\ge0$ for $i\ne j$, $\sum_{j=1}^n a^{ij}\ge0$ for $i\in I.$
	\end{enumerate}
	
	Since \eqref{I:eq:L-System} is invariant under adding constants to each
	component, we consider the solutions in the following Hilbert space:
	\[
	\mathcal H:=(\oH)^n,
	\qquad
	\oH:=\left\{v\in H^1(M):\int_Mv\,dV_g=0\right\},
	\]
	with the positive Hilbert norm
	$
	\|u\|_{\mathcal H}^2:=\sum_{i=1}^n
	\int_M|\nabla_gu_i|_g^2\,dV_g.
	$
	Weak solutions of the Liouville system \eqref{I:eq:L-System} are in one-to-one correspondence with critical points of the Euler--Lagrange functional: 
	\begin{equation*} J_\rho(u)=\frac12\sum_{i,j=1}^na^{ij}
		\int_M\nabla_gu_i\cdot\nabla_gu_j\,dV_g
		-\sum_{i=1}^n\rho_i\log\int_Mh_ie^{u_i}\,dV_g.
	\end{equation*}
	
	Liouville systems extend the scalar mean field equation and arise in
	conformal geometry, vortex theory, and models of interacting species;
	see \cite{CK1995,CSW1997,LinZhang2011}. The analysis of blow-up solutions
	depends on the finite mass entire system on $\mathbb R^2$.
	Chanillo and Kiessling \cite{CK1995} established an early symmetry theorem,
	and Chipot, Shafrir, and Wolansky \cite{CSW1997} developed the existence
	and symmetry theory. Lin and Zhang \cite{LinZhang2010} proved uniqueness
	for prescribed masses, up to translation and dilation, and obtained local
	estimates near an isolated blow-up point. Unlike the scalar equation,
	the entire system has a continuum of mass vectors. Comparing the entire
	profiles at different blow-up points is therefore part of the analysis.
	
	For $N\in\mathbb N$, set
	\begin{equation}\label{I:def:Lam_rho}
		\Lambda_{I,N}(\rho)
		:=4\sum_{i\in I}\frac{\rho_i}{2\pi N}
		-\sum_{i,j\in I}a_{ij}\frac{\rho_i\rho_j}{(2\pi N)^2},
		\qquad
		\Gamma_N:=\{\rho\in(\mathbb R_+)^n:\Lambda_{I,N}(\rho)=0\}.
	\end{equation}
	Lin and Zhang \cite{LinZhang2011} proved uniform bounds on compact subsets
	of the chambers
	\[
	O_N:=\left\{\rho\in(\mathbb R_+)^n:
	8\pi(N-1)\sum_i\rho_i
	<\sum_{i,j}a_{ij}\rho_i\rho_j
	<8\pi N\sum_i\rho_i\right\}
	\]
	and computed the corresponding Leray--Schauder degree.
	At a critical parameter $\rho^*\in\Gamma_N$, blow-up solutions have
	$N$ distinct concentration points. Lin and Zhang
	\cite{LinZhang2013} obtained sharp estimates for the blow-up location
	and parameter defect in the one bubble case. Huang and Zhang
	\cite{HuangZhang2022} treated multiple bubbles, including comparisons of
	their heights and local masses and the leading parameter expansion.
	Cheng, Li, and Zhang \cite{ChengLiZhangConstruction2025} constructed
	multiple bubble solutions on a flat torus by a Lyapunov--Schmidt reduction.
	Their subsequent work \cite{ChengLiZhangUniqueness2025} established local
	uniqueness and nondegeneracy under suitable assumptions, using refined
	profile estimates. For related constructions in the scalar case, see
	\cite{EspositoFigueroa2014}.
	
	The present paper concerns the approximation of  blow-up
	solutions in coordinates adapted to the variational problem at resonance.
	The local estimates in \cite{LinZhang2013,HuangZhang2022} compare solutions
	with entire profiles, while the construction in
	\cite{ChengLiZhangConstruction2025} uses center and scale parameters.
	For the method of critical points at infinity
	\cite{Bahri1989,BahriCoron1991}, we also use amplitude parameters and
	require the remainder to be orthogonal to every tangent direction of the
	approximate manifold. We obtain this decomposition by minimizing the
	distance in the product Dirichlet norm. The amplitude directions have
	norm of order $\sqrt{\log\lambda_t}$, and their contribution to the local
	mass contains $\lambda_t^{\sigma_i b_{i,t}}$. The estimate
	$L_k|b^k|\to0$, proved below, controls this contribution. These results
	provide the analytic foundation for the resonant degree problem.
	
	We fix $\rho^*\in\Gamma_N$ and write
	\[
	\sigma_i:=\frac{\rho_i^*}{2\pi N},\qquad
	m_i:=\sum_j a_{ij}\sigma_j,\qquad
	\widehat m:=\min_i m_i,\qquad
	\widehat I:=\{i:m_i=\widehat m\}.
	\]
	The normalized entire profile satisfies
	\[
	-\Delta v_i=\sum_j a_{ij}e^{v_j}\quad\text{in }\mathbb R^2,
	\qquad \frac1{2\pi}\int_{\mathbb R^2}e^{v_i}=\sigma_i,
	\qquad \max_i v_i(0)=0.
	\]
	Its common center is the origin. We have $m_i>2$, and either
	$2<\widehat m<4$ or $m_i=4$ for all $i$. Set
	\[
	q_i:=m_i-2,\qquad q:=(q_1,\ldots,q_n)^T,
	\]
	and
	\begin{equation*}
		\Psi_{\widehat m}(\lambda):=
		\begin{cases}
			\lambda^{2-\widehat m},&2<\widehat m<4,\\
			\lambda^{-2}\log\lambda,&\widehat m=4.
		\end{cases}
	\end{equation*}
	
	Let $G^g$ be the Green function of $-\Delta_g$ with mean zero. In a
	normalized isothermal chart $y_p$, we use the convention
	\[
	G^g(x,p)=-\frac1{2\pi}\chi_p(x)\log|y_p(x)|+H^g(x,p),
	\qquad R^g(p):=H^g(p,p),
	\]
	where $\chi_p=1$ near $p$. Set  $I':=\{1,\ldots,N\}$ and
	$\Delta:=\{p\in M^N:p_s=p_t\text{ for some }s\ne t\}$.
	For $p\in M^N\setminus\Delta$, define
	\begin{align}
		\widetilde G_t(x;p)&:=H^g(x,p_t)+\sum_{s\ne t}G^g(x,p_s),
		\label{I:def:tilde-G}\\
		\rf_{p,t}^i(x)&:=h_i(x)e^{2\pi m_i\widetilde G_t(x;p)}.
		\label{I:def:rf_i}
	\end{align}
	We abbreviate $\rf_p^i(p_t):=\rf_{p,t}^i(p_t)$; spatial derivatives
	always act on $x$ with the centers fixed. The reduced functional is
	\begin{equation}\label{I:def:reduce-F}
		\ff(p):=\sum_{t=1}^N\sum_{i=1}^n\rho_i^*
		\left[\log h_i(p_t)+\pi m_i
		\left(R^g(p_t)+\sum_{s\ne t}G^g(p_t,p_s)\right)\right].
	\end{equation}
	It satisfies
	\begin{equation}\label{I:eq:gradient-F}
		\nabla_{p_t}\ff(p)=\sum_i\rho_i^*\nabla_x\log\rf_{p,t}^i(p_t).
	\end{equation}
	
	Fix an orthonormal basis $e^1,\ldots,e^{N-1}$ of
	$E_0:=\{\delta\in\mathbb R^N:\sum_t\delta_t=0\}$, and write
	\begin{equation}\label{I:def:delta-star}
		\delta_t^*:=\sum_{\alpha=1}^{N-1}e_t^\alpha r_\alpha^*.
	\end{equation}
	\begin{definition}
		A pair $(p^*,r^*)$ is scale balanced if
		\begin{equation}\label{I:eq:balance-weights}
			\rf_{p^*,t}^i(p_t^*)e^{q_i\delta_t^*}
			=\rf_{p^*,s}^i(p_s^*)e^{q_i\delta_s^*}
			\qquad(i\in I,\ s,t\in I').
		\end{equation}
		Equivalently,
		\begin{equation*}
			\sum_{t=1}^N e_t^\alpha\log\rf_{p^*,t}^i(p_t^*)
			+q_i r_\alpha^*=0
			\qquad(i\in I,\ 1\le\alpha\le N-1).
		\end{equation*}
		Since $q_i>0$, $r^*$ is unique whenever it exists.
	\end{definition}
	Set
	\begin{equation}\label{I:def:balanced-critical-set}
		\mathcal K_{\mathrm{bal}}
		:=\{(p^*,r^*):\nabla\ff(p^*)=0,\ (p^*,r^*)\text{ is scale balanced}\}.
	\end{equation}
	Permutations act simultaneously on $p^*$ and $\delta^*$.
	For $N=1$, the relative scale conditions are empty.
	We use the following assumptions:
	\begin{enumerate}[label=(A\arabic*),ref=(A\arabic*)]
		\item $\ff$ is Morse on $M^N\setminus\Delta$;
		\item $\mathcal L_\infty(p^*,r^*)\ne0$ for every
		$(p^*,r^*)\in\mathcal K_{\mathrm{bal}}$.
	\end{enumerate}
	Here $\mathcal L_\infty$ is the leading coefficient defined in
	\eqref{I:def:L-balanced}. Assumption (A1) is used for matching and
	finiteness; (A2) is needed for compactness at $\rho^*$.
	
	We now introduce the approximate solutions. Set
	\[
	W_{i,p_t,\lambda_t}(x):=v_i(\lambda_t y_{p_t}(x))+2\log\lambda_t
	\]
	in the isothermal chart $(y_{p_t}, \mathcal U(p_t))$ at $p_t$. Its projection $PW_{i,p_t,\lambda_t}$ is determined by
	\[
	-\Delta_g PW_{i,p_t,\lambda_t}
	=\chi_{p_t}e^{-\varphi_{p_t}}e^{W_{i,p_t,\lambda_t}}
	-\int_M\chi_{p_t}e^{-\varphi_{p_t}}e^{W_{i,p_t,\lambda_t}}\,dV_g,
	\qquad \int_M PW_{i,p_t,\lambda_t}\,dV_g=0,
	\]
	where $g=e^{\varphi_{p_t}}|dy_{p_t}|^2$ in the chart and $\chi_{p_t}$
	is a smooth cutoff function supported in $\mathcal U(p_t)$ and equal to $1$
	near $p_t$.  For
	$b=(b_{i,t})\in\mathbb R^{nN}$, let
	\[
	U_{p,\lambda,b}:=(U_{1,p,\lambda,b},\ldots,U_{n,p,\lambda,b}),
	\qquad
	U_{i,p,\lambda,b}:=\sum_{t=1}^N
	\left(\sum_j a_{ij}PW_{j,p_t,\lambda_t}
	+b_{i,t}PW_{i,p_t,\lambda_t}\right).
	\]
	These profiles form the approximate manifold defined in
	\eqref{def:proj-bubble}. Its normal space is
	\[
	\begin{aligned}
		E^N_{p,\lambda,b}
		:=\Big\{w\in\mathcal H: &
		\langle w,\partial_{b_{i,t}}U\rangle_{\mathcal H}
		= 0, \forall  i\in I, t\in I'; \langle w,\partial_{\lambda_t}U\rangle_{\mathcal H}
		=0, \forall t\in I';\\
		&	\langle w,\partial_{(p_t)_l}U\rangle_{\mathcal H}=0, \forall   t\in I', l\in \{1,2\}
		\Big \} .
	\end{aligned}
	\]
	We use the logarithmic scale parametrization
	\begin{equation*}
		\lambda(\ell):=(e^{\ell_1},\ldots,e^{\ell_N}),
		\qquad \ell\in\mathbb R^N.
	\end{equation*}
	
	\begin{definition}\label{I:def:admissible-path}
		A parameter sequence $\rho^k\to\rho^*$ is called admissible if
		\[
		\rho^k=\rho^*+s_k\dot\rho+o(|s_k|),\qquad s_k\to0,
		\qquad \dot\rho_i\ne0\text{ with the same sign for all }i.
		\]
	\end{definition}
	It follows that 
	\[
	D\Lambda_{I,N}(\rho^*)[\dot\rho]
	=-\frac1{\pi N}\sum_iq_i\dot\rho_i\ne0,
	\qquad |\Lambda_{I,N}(\rho^k)|\asymp|\rho^k-\rho^*|.
	\]
	The lower side is $\Lambda_{I,N}>0$, adjacent to $O_N$.
	
	Our first result gives the quantitative approximation and its unique
	orthogonal decomposition.
	\begin{theorem}\label{I:thm:intro-coarse-entry}
		Assume (A1). Let $u^k\in\mathcal H$ be blow-up solutions of
		\eqref{I:eq:L-System} with $\rho^k\to\rho^*\in\Gamma_N$.
		Suppose either $\rho^k\equiv\rho^*$ or $\rho^k$ is admissible with
		$\Lambda_{I,N}(\rho^k)>0$. For $N=1$, either admissible side is allowed.
		After passing to a subsequence and relabeling the bubbles, there exist
		parameters $(\widehat p^k,\widehat\lambda^k)$ with comparable scales such that
		\[
		\|u^k-U_{\widehat p^k,\widehat\lambda^k,0}\|_{W^{2,p}(M)^n}
		\le Ce^{-\kappa_p\widehat L_k},\qquad
		\widehat L_k:=\frac1N\sum_t\log\widehat\lambda_t^k,
		\qquad \kappa_p:=\frac{\beta-2(p-1)}p>0,
		\]
		for sufficiently small fixed $\beta>0$ and $1<p<1+\beta/2$.
		
		For $\widehat\ell^k_t=\log \widehat \lam_t^k, \widehat z^k=(\widehat p^k,\widehat\ell^k,0)$, the problem
		\[
		\min_{(p,\ell,b)\in\overline{\mathcal D_{\widehat z^k}(r_*)}}
		\|u^k-U_{p,\lambda(\ell),b}\|_{\mathcal H}^2
		\]
		has a unique interior minimizer for each ordering. Here
		$\mathcal D_{\widehat z^k}(r_*)$ is the normalized parameter neighborhood
		in \eqref{I:def:modulation-box}. The resulting decomposition
		\[
		u^k=U_{p^k,\lambda^k,b^k}+w^k,\qquad
		w^k\in E^N_{p^k,\lambda^k,b^k},
		\]
		is unique up to permutation. With $L_k:=N^{-1}\sum_t\log\lambda_t^k$,
		\[
		\min_t\lambda_t^k\to\infty,\qquad
		L_k|b^k|+\|w^k\|_{\mathcal H}+\max_{i,t}|\tau_{i,t}^k|\to0,
		\qquad
		p^k\to p^*,\quad \log\lambda_t^k-L_k\to\delta_t^*,
		\]
		where $(p^*,r^*)\in\mathcal K_{\mathrm{bal}}$ and
		$\tau_{i,t}^k$ is defined in \eqref{I:def:tau_ijt}.
	\end{theorem}
	
	Retaining the entire mass vector at each concentration point gives an
	exact mass variance identity. Together with matching on fixed circles
	and integrated local estimates, it yields the parameter expansion.
	\begin{theorem}
		For the admissible sequences in Theorem~\ref{I:thm:intro-coarse-entry},
		the minimizing scales satisfy, for every $t\in I'$,
		\[
		\Lambda_{I,N}(\rho^k)
		=-\mathcal L_t(p^*)\Psi_{\widehat m}(\lambda_t^k)
		+o\bigl(\Psi_{\widehat m}(\lambda_t^k)\bigr).
		\]
		Equivalently,
		\[
		\Lambda_{I,N}(\rho^k)
		=-\mathcal L_\infty(p^*,r^*)\Psi_{\widehat m}(e^{L_k})
		+o\bigl(\Psi_{\widehat m}(e^{L_k})\bigr).
		\]
		The coefficients $\mathcal L_t(p^*)$ and $\mathcal L_\infty(p^*,r^*)$ are given in \eqref{I:def:tail-coefficients} and
		\eqref{I:def:L-balanced}, respectively. 
	\end{theorem}
	Thus, when  $\mathcal L_\infty(p^*,r^*)\ne0$, the sign of
	$\Lambda_{I,N}(\rho^k)$ is opposite to that of
	$\mathcal L_\infty(p^*,r^*)$ for large $k$.
	
	The approximation is uniform along each fixed admissible path.
	\begin{theorem}
		Assume (A1), and fix a $C^1$ path $\rho(s)$ with
		$\rho(0)=\rho^*$ and all components of $\rho'(0)$ nonzero and of the
		same sign. There exist $\varepsilon,R>0$ and finitely many pairwise
		disjoint bubbling neighborhoods, invariant under permutations, such that
		every solution $u$ of \eqref{I:eq:L-System} with parameter $\rho(s)$ satisfying
		\[
		0<|s|<\varepsilon,\qquad
		\Lambda_{I,N}(\rho(s))>0,\qquad \|u\|_{L^\infty}>R
		\]
		belongs to exactly one of these neighborhoods. Its decomposition in
		Theorem~\ref{I:thm:intro-coarse-entry} is unique up to permutation.
	\end{theorem}
	
	At the fixed critical parameter, the same mass identity gives compactness.
	\begin{corollary}
		Assume (A1)--(A2). The solution set of \eqref{I:eq:L-System} at $\rho^*$
		is bounded in $L^\infty(M)^n$ and compact in $C^{2,\alpha}(M)^n$ for
		every $0<\alpha<1$.
	\end{corollary}

	\subsection*{Organization of the paper}
	Subsection~\ref{I:sec:bubbling-chart} introduces the projected profiles
	and the minimizing decomposition. Subsection~\ref{I:sec:coarse-blowup}
	fixes the normalization, and Proposition~\ref{I:prop:refined-local-Harnack}
	records the refined local estimates.
	Section~\ref{I:sec:sharp-defect} proves the mass variance identity,
	matching estimates, and sharp parameter expansion, ending with the
	summary in Subsection~\ref{I:sec:multiple-coupling}.
	Section~\ref{I:sec:entry-theorem} estimates the local masses and completes
	the approximation theorem. Section~\ref{I:sec:uniform-capture} proves
	uniform separation, finiteness, and compactness.
	
	\section{Preliminaries}
	\subsection{Entire profiles and projected bubbling coordinates}\label{I:sec:bubbling-chart}
	
	For every nonempty $J\subset I$, set
	\begin{equation}\label{I:def:Lam_J}
		\Lambda_J(\sigma):=4\sum_{i\in J}\sigma_i
		-\sum_{i,j\in J}a_{ij}\sigma_i\sigma_j.
	\end{equation}
	The constant terms at infinity of the normalized profile are
	\begin{equation}\label{I:def:J-i}
		\fJ_i:=v_i(0)+\sum_{j\in I}\frac{a_{ij}}{2\pi}
		\int_{\mathbb R^2}e^{v_j(y)}\log|y|\,dy.
	\end{equation}
	The limit system of \eqref{I:eq:L-System} is as follows:
	\begin{equation}
		\label{I:eq:entire} \left\{\begin{array}{ll}
			-\Delta v_i=\sum_{j=1}^n a_{i j} e^{v_j}   & \text { in  }  \mathbb{R}^2, \\
			\frac{1}{2 \pi} \int_{\mathbb{R}^2} e^{v_i}\, dy =\sigma_i & \text { for } i \in I.
		\end{array}\right.
	\end{equation}
	Under \ref{I:assump:H1}, every entire solution of
	\eqref{I:eq:entire} is radially symmetric with respect to a common point
	by Chipot--Shafrir--Wolansky \cite{CSW1997}; the exact statement is
	also recalled as Theorem~A in \cite{LinZhang2010}. The earlier symmetry
	theorem of \cite{CK1995} has additional diagonal hypotheses and is
	not used to cover matrices with zero diagonal entries.  Its mass vector is the fixed vector
	$\sigma=(\sigma_1,\ldots,\sigma_n)$ introduced above and satisfies
	\[
	\Lambda_I(\sigma)=0,
	\qquad
	\Lambda_J(\sigma)>0
	\quad\text{for every }\varnothing\ne J\subsetneq I,
	\]
	where $\Lambda_J$ is defined in \eqref{I:def:Lam_J}.
	
	Moreover, \cite[Lemma 3.1]{HuangZhang2022} implies that every radial solution $v=\left(v_1, \ldots, v_n\right)$ of \eqref{I:eq:entire} with $\max_{i\in I} v_i(0)=0$ satisfies:
	\begin{equation}
		\label{I:eq:est-entire} 	\left\{\begin{aligned}
			v_i(y)&=-m_i\log |y|+\fJ_i-\sum_{j\in I} \frac{a_{i j} e^{\fJ_j}}{(m_j-2)^2}|y|^{2-m_j}+\cO(|y|^{2-\widehat m-\delta_{\mathrm{as}}}), \\
			\nabla v_i(y)&=-m_i  \frac{y}{|y|^2}+ \sum_{j\in I} \frac{a_{ij} e^{\fJ_j}}{m_j-2}\frac{y}{|y|^{m_j}}+\cO(|y|^{1-\widehat m-\delta_{\mathrm{as}}}),\quad \forall i\in I,
		\end{aligned}\right.
	\end{equation}
	for $|y|$ large, as $|y|\to \infty$.  Here $m_i$ and $\widehat m$ are the fixed resonant quantities defined above, $\fJ_i$ is given by \eqref{I:def:J-i}, and $\delta_{\mathrm{as}}>0$ is the fixed decay margin;
	\begin{equation*}
		\left\{\begin{array}{l}
			v_i(y)=v_i(0)-\frac{1}{4} \sum_{j=1}^n a_{i j} e^{v_j(0)}|y|^2+ \cO\left(|y|^4\right),\\
			\nabla v_i(y)=-\frac{1}{2} \sum_{j=1}^n a_{i j} e^{v_j(0)} y+\cO\left(|y|^3\right),
		\end{array}\right.
	\end{equation*}
	as $|y|\to0$.  Throughout this paper, $\delta_{\mathrm{as}}>0$ denotes the fixed decay margin in \eqref{I:eq:est-entire}; the symbol $\delta$ is reserved later for the coarse exact chart.
	%--------------
	\begin{lemma}
		\label{I:lem:radial-family-derivatives}
		Let
		\[
		\Pi:=\left\{
		a\in(0,\infty)^n:\Lambda_I(a)=0
		\right\},
		\]
		and fix $\sigma\in\Pi$. There exists a relatively open neighborhood
		$U$ of $\sigma$ in $\Pi$, with compact closure $K\subset\Pi$, such
		that the normalized radial family
		\[
		-\Delta V_i(a)=\sum_j a_{ij}e^{V_j(a)}
		\quad\text{in }\mathbb R^2,\qquad
		\frac1{2\pi}\int_{\mathbb R^2}e^{V_i(a)}=a_i,
		\qquad V_n(a;0)=v_n(0),
		\]
		satisfies, uniformly for $a\in K$,
		\[
		e^{V_i(a;y)}\le C(1+|y|)^{-2-\eta^{\#}},
		\qquad
		|D_a^jV_i(a;y)|
		\le C_j\bigl(1+\log(2+|y|)\bigr),
		\quad j=1,2,
		\]
		for some $\eta^{\#}>0$. Here $D_a$ denotes differentiation in the
		local coordinates $(a_1,\ldots,a_{n-1})$ on $\Pi$.
		
		The asymptotic constants $c_i(a)$ are $C^2$. Moreover, for $j=0,1,2$,
		\[
		\left|
		D_a^j\left(\frac{e^{V_i(a;y)}}{2\pi a_i}\right)
		\right|
		\le C(1+|y|)^{-2-\eta^{\#}}
		\bigl(1+\log(2+|y|)\bigr)^2.
		\]
		This majorant remains integrable after multiplication by
		$(1+|\log|y||)^q$ for every fixed $q\ge0$.
	\end{lemma}
	
	\begin{proof}
		We first use the initial value coordinates
		\[
		V_i(\beta;0)=\beta_i\quad(i<n),
		\qquad V_n(\beta;0)=v_n(0),
		\]
		near $\beta^*=(v_1(0),\ldots,v_{n-1}(0))$.
		Existence and uniqueness of the normalized profiles follow from
		\cite[Theorem~1.1 and Theorem~B]{LinZhang2010}.
		
		Set
		\[
		m_i^r(\beta):=\sum_j a_{ij}
		\int_0^r s e^{V_j(\beta;s)}\,ds.
		\]
		Since $m_i^r(\beta^*)\to(A\sigma)_i>2$, there exist $R\ge2$ and
		$\eta^{\#}>0$ such that
		$m_i^R(\beta)>2+\eta^{\#}$ for $\beta$ near $\beta^*$.
		Here we use smooth ODE dependence on bounded intervals.
		The identities
		\[
		rV_i'(\beta;r)=-m_i^r(\beta),
		\qquad
		m_i^r(\beta)\ge m_i^R(\beta)\quad(r\ge R)
		\]
		then give
		\[
		e^{V_i(\beta;r)}\le C(1+r)^{-2-\eta^{\#}}.
		\]
		All constants are uniform on a smaller compact initial value
		neighborhood.
		
		For $Y=\partial_{\beta_l}V$, differentiation on finite intervals gives
		\[
		Y_i(r)=Y_i(R)+RY_i'(R)\log\frac rR
		-\sum_j a_{ij}\int_R^r
		s e^{V_j(s)}Y_j(s)\log\frac rs\,ds.
		\]
		Setting
		\[
		\|Y\|_{R,T}:=
		\max_i\sup_{R\le r\le T}
		\frac{|Y_i(r)|}{1+\log r},
		\]
		we obtain
		\[
		\|Y\|_{R,T}
		\le C_R+
		C\int_R^\infty
		s(1+s)^{-2-\eta^{\#}}(1+\log s)\,ds\,
		\|Y\|_{R,T}.
		\]
		Increase $R$ so that the coefficient is at most $1/2$ and let
		$T\to\infty$. Thus
		$|D_\beta V_i(r)|\le C(1+\log r)$.
		This is the logarithmic growth estimate used in
		\cite[Lemma~2.1 and Section~3]{LinZhang2010}, with uniform constants.
		
		For $Z=\partial_{\beta_l}\partial_{\beta_m}V$, the same equation
		has the additional source
		$e^{V_j}\partial_{\beta_l}V_j\partial_{\beta_m}V_j$.
		Since
		\[
		\int_R^\infty
		s(1+s)^{-2-\eta^{\#}}(1+\log s)^2\,ds<\infty,
		\]
		the same absorption yields
		$|D_\beta^2V_i(r)|\le C(1+\log r)$.
		Consequently,
		\[
		|D_\beta^j e^{V_i(\beta;y)}|
		\le C(1+|y|)^{-2-\eta^{\#}}
		\bigl(1+\log(2+|y|)\bigr)^2,
		\qquad j=0,1,2.
		\]
		
		Dominated differentiation shows that the mass map
		\[
		\beta\longmapsto
		\left(
		\int_0^\infty r e^{V_1(\beta;r)}\,dr,\ldots,
		\int_0^\infty r e^{V_{n-1}(\beta;r)}\,dr
		\right)
		\]
		is $C^2$. Its Jacobian is nonsingular by the proof of
		\cite[Theorem~1.1, Section~3]{LinZhang2010}.
		The inverse function theorem therefore gives a $C^2$ inverse.
		Choose $U$ with compact closure inside this mass chart.
		The chain rule gives the stated mass derivative estimates on $K$.
		Since $\inf_{a\in K}a_i>0$, they also give the bounds for the
		normalized densities.
		
		Finally, dominated differentiation in
		\[
		c_i(a)=V_i(a;0)+\sum_j a_{ij}
		\int_0^\infty r e^{V_j(a;r)}\log r\,dr
		\]
		proves that $c_i\in C^2(K)$. The same majorants justify
		differentiation of every fixed logarithmic moment.
	\end{proof} 
	
	%--------------
	\begin{lemma}
		Let $V(a)$ be the normalized radial family of
		Lemma~\ref{I:lem:radial-family-derivatives}, with
		$V_n(a;0)=v_n(0)$. Set 
		\[
		m(a):=Aa,\qquad q_i(a):=m_i(a)-2,\qquad q:=q(\sigma).
		\]
		Fix the index sets at the base mass $\sigma$:
		\[
		\begin{aligned}
			\mathscr K
			&:=\left\{\mathbf k\in\mathbb N_0^n:
			|\mathbf k|\ge1,\ \mathbf k\cdot q\le1\right\},\\
			\mathscr K_i
			&:=\left\{\mathbf k\in\mathbb N_0^n:
			q_i+\mathbf k\cdot q\le1\right\}.
		\end{aligned}
		\]
		After shrinking the compact mass neighborhood $K$ of $\sigma$,
		there exist $\eta>0$ and $C^2$ coefficients
		$c_{i,\mathbf k}(a)$ and $d_{i,\mathbf k}(a)$ such that
		\begin{equation*}
			%\label{I:eq:finite-tail-ladder}
			V_i(a;r)
			=-m_i(a)\log r+\fJ_i(a)
			+\sum_{\mathbf k\in\mathscr K}
			c_{i,\mathbf k}(a)r^{-\mathbf k\cdot q(a)}
			+\mathcal R_i(a;r)
		\end{equation*}
		and
		\begin{equation}\label{I:eq:finite-mass-tail-ladder}
			\int_{|y|>r}e^{V_i(a;y)}\,dy
			=\sum_{\mathbf k\in\mathscr K_i}
			d_{i,\mathbf k}(a)
			r^{-(q_i(a)+\mathbf k\cdot q(a))}
			+\mathcal S_i(a;r).
		\end{equation}
		Uniformly for $a\in K$ and $r\ge2$, it holds 
		\[
		|D_a^j(r\partial_r)^l\mathcal R_i(a;r)|
		+|D_a^j(r\partial_r)^l\mathcal S_i(a;r)|
		\le Cr^{-1-\eta},
		\qquad 0\le j,l\le2.
		\]
		Here $D_a$ denotes differentiation in local mass coordinates.
		If $q_i(\sigma)\le1$, then
		\begin{equation}
			\label{eq:d_i0}  d_{i,0}(a)=\frac{2\pi e^{\fJ_i(a)}}{q_i(a)}.
		\end{equation}
	\end{lemma}
	
	\begin{proof}
		The leading expansion follows from
		\cite[Lemma~2.2]{LinZhang2013}, and its first correction is given in
		\cite[Lemma~3.1]{HuangZhang2022}. We iterate the same radial argument.
		
		Set $w_i(a;r):=V_i(a;r)+m_i(a)\log r-\fJ_i(a)$.
		The radial equation gives
		\begin{equation}\label{I:eq:radial-volterra}
			w_i(a;r)
			=-\sum_j a_{ij}e^{\fJ_j(a)}
			\int_r^\infty\frac{ds}{s}
			\int_s^\infty t^{-1-q_j(a)}e^{w_j(a;t)}\,dt.
		\end{equation}
		Successive Taylor expansion of $e^{w_j}$ produces powers
		$r^{-\mathbf k\cdot q(a)}$, since
		\[
		\int_r^\infty\frac{ds}{s}
		\int_s^\infty t^{-1-q_j(a)-\gamma}\,dt
		=\frac{r^{-q_j(a)-\gamma}}{(q_j(a)+\gamma)^2},
		\qquad \gamma\ge0.
		\]
		All denominators are uniformly positive. Each iteration improves
		the remainder exponent by at least
		$\min_{a\in K,j}q_j(a)>0$, so finitely many iterations suffice.
		
		The exponents at $\sigma$ form locally finite sets. Choose $\eta>0$
		and shrink $K$ so that every omitted profile or mass exponent exceeds
		$1+4\eta$. The index sets remain fixed as $a$ varies.
		Lemma~\ref{I:lem:radial-family-derivatives} justifies two mass
		derivatives in \eqref{I:eq:radial-volterra}. These derivatives introduce
		at most $(\log r)^2$, while $r\partial_r$ preserves each power.
		The same iteration gives the differentiated remainder estimates,
		using
		\[
		r^{-1-4\eta}(1+\log r)^2\le Cr^{-1-\eta}.
		\]
		Finally, expand
		$e^{V_i}=e^{\fJ_i}r^{-2-q_i}e^{w_i}$
		and integrate in polar coordinates. This gives
		\eqref{I:eq:finite-mass-tail-ladder} and \eqref{eq:d_i0}. 
	\end{proof}

	%----------------------------------
	We next construct the approximate solutions.  Let
	\begin{equation*}
		\chi(s)=\begin{cases}
			1 & |s|\leq \frac 1 2, \\
			0 & |s|\geq 1,
		\end{cases}
	\end{equation*}
	be a smooth cut-off function in $C^\infty ( \R, [0,1])$.
	
	For any $p\in M$, we define the Green function $G^g(\cdot,p)=G^g_p$ by the following equations:
	\begin{equation*}
		\begin{cases}
			-\Delta_g G^g(x,p)= \delta_p -1, & x\in M, \\
			\int_M G^g(x,p)\, dV_g(x)=0, &
		\end{cases}
	\end{equation*}
	where $\delta_p$ is the Dirac mass on $M$ concentrated at $p$.
	
	By the existence of isothermal coordinates for Riemann surfaces \cite{Chern1955}, and
	following \cite[Section~2]{EspositoFigueroa2014}, we choose local
	families of charts
	\[
	y_p:U(p)\longrightarrow B_{2r_0}(0),
	\qquad y_p(p)=0,
	\]
	with a common radius $r_0>0$, such that
	\[
	g=e^{\widehat\varphi_p(y)}(dy_1^2+dy_2^2),
	\qquad
	\widehat\varphi_p(0)=0,
	\qquad
	\nabla\widehat\varphi_p(0)=0.
	\]
	The charts and conformal factors depend smoothly on $p$ locally. 
	On a finite cover, their derivatives of each fixed order are uniformly
	bounded on smaller coordinate disks. Moreover,
	\[
	-\Delta\widehat\varphi_p(y)
	=2K_g(y_p^{-1}(y))e^{\widehat\varphi_p(y)}, 
	\]
	where $K_g$ is the Gaussian curvature of $(M,g).$
	Set
	\[
	\varphi_p:=\widehat\varphi_p\circ y_p,
	\qquad
	U_r(p):=y_p^{-1}(B_r(0)),
	\quad 0<r<2r_0.
	\]
	
	Choose $\chi\in C_c^\infty([0,1))$ with $0\le\chi\le1$ and
	$\chi\equiv1$ on $[0,1/2]$, and define
	\[
	\chi_p(x):=
	\begin{cases}
		\chi(|y_p(x)|/r_0),&x\in U(p),\\
		0,&x\notin U(p).
	\end{cases}
	\]
	Then $\chi_p\equiv1$ on $U_{r_0/2}(p)$ and
	$\operatorname{supp}\chi_p\subset U_{r_0}(p)$.
	With this convention, the regular part of the Green function is
	\[
	H^g(x,p)
	:=G^g(x,p)-\frac1{2\pi}\chi_p(x)
	\log\frac1{|y_p(x)|}.
	\]
	
	Then, $H^g_p$ is a smooth solution of the following equations:
	\begin{equation*}
		\left\{\begin{aligned}
			-\Delta_g H^g_p
			&=
			-\frac 1 {2\pi} (\Delta_g \chi_{p}) \log  {|y_p|}
			-\frac 1 {\pi}\la\nabla\chi_{p}, \nabla\log {|y_p|}\ra_g 	-1 & \text{ in } M \\
			\int_M H^g_p \,dV_g&=  \frac 1 {2\pi}\int_{M}  \chi_{p}\log {|y_p|}\,  dV_g &&
		\end{aligned}\right..\end{equation*}
	Elliptic regularity gives uniform smooth bounds for $H_p^g$.
	The Robin function
	$
	R^g(p):=H^g(p,p)
	$
	is independent of the cutoff and the normalized isothermal chart.
	Indeed, $\chi_p$ is equal to $1$ near $p$, and any two normalized charts
	$y_p$ and $\widetilde y_p$ satisfy
	$
	\lim_{x\to p}
	\log\frac{|\widetilde y_p(x)|}{|y_p(x)|}=0.
	$

	For a critical point $p^*\in M^N\setminus\Delta$ of $\ff$, define
	\[ M_{p^*,\tau}:=\Big\{ p\in M^N \setminus \Delta:  \max_{t\in I'} d_g(p^*_t, p_t)\leq \tau  \text{ and }  U(p_t)\cap U(p_s)=\emptyset, \forall s\neq t\Big\},\]
	where $d_g(\cdot,\cdot)$ denotes the geodesic distance on $M$ induced by $g$.
	For sufficiently small $\tau>0$, the isothermal charts $(U(p_t),y_{p_t})$ may be chosen so that $y_{p_t}$ and $\widehat{\varphi}_{p_t}$ depend smoothly on $p$; moreover, one may use a uniform radius $r_{p_t}=r_0$ for every $t\in I'$ and $p\in M_{p^*,\tau}$.
	
	Given $\rho^*\in\Gamma_N$ and $p:=(p_1,\ldots,p_N)\in M_{p^*,\tau}$, for $i\in I$ and $t\in I':=\{1,\ldots,N\}$, we construct the bubbles as follows:
	\[ W_{i, p_t,\lam_t}(x):= \begin{cases}
		v_i (  \lam_{t}y_{p_t}(x))
		+ 2 \log \lam_{t},  \quad &  x\in U(p_t), \\
		0, \quad & x\in M\setminus U(p_t),
	\end{cases}  \]
	where $v=(v_1,\ldots,v_n)$ is the normalized radial entire solution of \eqref{I:eq:entire} associated with the fixed mass vector $\sigma=(\sigma_1,\ldots,\sigma_n)$, so that $\max_{i\in I}v_i(0)=0$.
	
	We project $W_{i,p_t,\lam_t}$ onto $\oH$ through
	\begin{equation*}
		\left\{ 	\begin{aligned}
			-\Delta_g  PW_{i, p_t,\lam_t} &=  \chi_{p_t} e^{-\varphi_{p_t}} e^{W_{i, p_t,\lam_t}}    - \overline{ \chi_{p_t} e^{-\varphi_{p_t}} e^{W_{i, p_t,\lam_t}}  }  && \text{ in } M,\\
			\int_M   PW_{i, p_t,\lam_t} \, dV_g &=0,&&
		\end{aligned} \right.
	\end{equation*}
	where the average is denoted by $\bar f:=\int_M f\,dV_g$ for $f\in L^1(M)$.
	We also define
	\begin{equation}
		\label{I:eq:proj-B2}
		\left\{ 	\begin{aligned}
			-\Delta_g PW^*_{i, p_t,\lam_t} &= \sum_{j\in I} a_{ij} \chi_{p_t} e^{-\varphi_{p_t}} e^{W_{j, p_t,\lam_t}}    - \overline{ \sum_{j\in I} a_{ij} \chi_{p_t} e^{-\varphi_{p_t}} e^{W_{j, p_t,\lam_t}}  }  && \text{ in } M,\\
			\int_M   PW^*_{i, p_t,\lam_t} \, dV_g &=0,&&
		\end{aligned} \right.
	\end{equation}
	By uniqueness for \eqref{I:eq:proj-B2},
	$PW^*_{i, p_t,\lam_t}=\sum_{j\in I} a_{ij}PW_{j,p_t,\lam_t}. $
	
	For $i=1,\ldots,n$, we define
	\begin{equation}
		\label{def:proj-bubble}  U_{i,p,\lam,b}
		:= \sum_{j=1}^n \sum_{t=1}^N \alpha_{i,j,t}\, PW_{j,p_t,\lam_t}=\sum_{t=1}^N  PW^*_{i, p_t,\lam_t}+ \sum_{t=1}^N b_{i,t} PW_{i,p_t,\lam_t},
	\end{equation}
	where the parameters  $p=(p_1,\ldots, p_N)\in M^N\setminus\Delta$, $\lam:=(\lam_1,\ldots,\lambda_N) \in (\R_+)^N ,$ and $\alpha_{i,j,t}:=a_{ij}+\delta_{ij}b_{i,t}.$
	Here $b=(b_{i,t})\in\R^{nN}$ and $\delta_{ij}$ is the Kronecker delta. We also define
	\begin{equation*}
		\alpha^{i,j,t}:=\sum_{k=1}^n a^{jk}\alpha_{i,k,t}
		=\delta_{ij}+a^{ji}b_{i,t}.
	\end{equation*}
	These coefficients satisfy the  following identity:
	\begin{equation*}
		\sum_{j\in I}\alpha^{i,j,t}m_j-2
		=q_i+\sigma_i b_{i,t}.
	\end{equation*}
	We also define 
	\begin{equation*}
		\rg_{p,\alpha}^i(x):=
		\exp\!\left(
		\sum_{j\in I}\sum_{t\in I'}2\pi m_j
		(\alpha^{i,j,t}-\delta_{ij})H^g(x,p_t)
		\right).
	\end{equation*}
	Thus $\rg_{p,\alpha}^i\equiv1$ when $b=0$.
	Let  \[U_{p,\lam,b}:=(U_{1,p,\lam,b},\ldots,U_{n, p,\lam,b})\in \cH.\]
	We consider $U_{p,\lam,b}$ as the approximate manifold and its tangent space is generated by
	\begin{align*} 
		X^b_{i,t}&:= \Big( 0, \ldots,0, \underbrace{  PW_{i,p_t,\lam_t}}_{i\text{-th entry}}, 0, \ldots, 0\Big)^\top, \\
		X^{\lam}_{t}&:= \Big( \sum_{j\in I} \alpha_{1,j,t} \partial_{\lam_t} PW_{j,p_t,\lam_t}, \ldots,  \sum_{j\in I} \alpha_{n,j,t} \partial_{\lam_t} PW_{j,p_t,\lam_t}\Big)^\top,  \\
		X^{p}_{t,l}&:=\Big( \sum_{j\in I} \alpha_{1,j,t} \partial_{(p_t)_l} PW_{j,p_t,\lam_t}, \ldots,  \sum_{j\in I} \alpha_{n,j,t} \partial_{(p_t)_l} PW_{j,p_t,\lam_t}\Big)^\top,  
	\end{align*}
	for $i\in I$, $t\in I'$, and $l=1,2$.
	
	We denote by
	$
	\mathbf 0:=(0)_{i\in I,\;t\in I'}
	$
	the zero vector in $\R^{nN}$.
	
	Given $N\in\mathbb N$, $\rho^*\in\Gamma_N$, and $p^*\in M^N\setminus\Delta$, we define the following ordered bubbling neighborhood: for $\eps>0$
	\begin{align*}
		V^\eta_{\rho^*}(N,\eps):= \Big\{u\in \cH: &
		\|\nabla J_{\rho^*}(u)\|<\eps;\ 
		\exists \lam\in(\R_+)^N, \, \max_{t\in I'}\lam_t^{-1}<\eps,\, 
		\max_{s,t\in I'}\frac{\lam_t}{\lam_s}<C,\\
		&\exists p=(p_1,\ldots,p_N),\,  p\in M_{p^*,\tau/2},\, 
		d_g(p_s,p_t)>2\eta\ \text{for all }s\ne t,\\
		& \text{s.t. }\max_{i\in I}\|u_i-U_{i,p,\lam,\mathbf0}\|<\eps\Big\},
	\end{align*}
	where $\tau>0, C>0$ and $\eta>0$ are fixed constants.  The base configuration $p^*$ is suppressed in the notation
	$V^\eta_{\rho^*}(N,\eps)$; globally we use the union of these charts over
	the finitely many separated critical configurations.

	Following \cite{BahriCoron1991}, we select the parameters by minimizing
	the squared distance in the product Dirichlet norm over an enlarged
	{\it ordered} bubbling chart. The minimizer exists and is unique within this chart. 
	
	Let
	\[
	\bar z=(\bar p,\bar\ell,0),\qquad
	\bar\ell_t:=\log\bar\lambda_t,\qquad
	\bar L:=\frac1N\sum_{t=1}^N\bar\ell_t,
	\]
	be preliminary ordered parameters in the definition of
	$V^\eta_{\rho^*}(N,\eps)$. For fixed small $r_*>0$, set
	\begin{equation}\label{I:def:modulation-box}
		\begin{aligned}
			\mathcal D_{\bar z}(r_*):=
			\Big\{
			(p,\ell,b):\;&p\in M_{p^*,\tau},\,\max_{t\in I'}\bar\lambda_t d_g(p_t,\bar p_t)<r_*,\, \max_{t\in I'}|\ell_t-\bar\ell_t|<r_*,\\
			& \sqrt{\bar L}\max_{i\in I,\,t\in I'}|b_{i,t}|<r_*
			\Big\}.
		\end{aligned}
	\end{equation}
	Setting  $\lambda(\ell):=(e^{\ell_1},\ldots,e^{\ell_N})$, we consider the following minimization problem:
	\begin{equation}\label{I:eq:exact-parameter-minimization}
		\min_{(p,\ell,b)\in\overline{\mathcal D_{\bar z}(r_*)}}
		\|u-U_{p,\lambda(\ell),b}\|_{\mathcal H}^2.
	\end{equation}

	%-------------------------------------
	\begin{lemma}\label{I:lem:unique}
		There exist $r_*,\eps_0>0$, independent of the preliminary ordered
		parameters, such that for every $0<\eps<\eps_0$ and
		$u\in V_{\rho^*}^\eta(N,\eps)$, problem
		\eqref{I:eq:exact-parameter-minimization} has a unique minimizer
		$(p,\ell,b)\in\mathcal D_{\bar z}(r_*)$. Any two preliminary representations yield the same minimizing parameters
		after a permutation of the bubble labels.
		
		Setting  $\lambda=\lambda(\ell)$, we have
		\[
		u=U_{p,\lambda,b}+w,\qquad
		w\in E^N_{p,\lambda,b},\qquad
		\|w\|_{\mathcal H}
		+\max_{i,t}\sqrt{\log\lambda_t}\,|b_{i,t}|
		\le C\eps.
		\]
		In particular, as $\eps\to0$,
		\begin{equation}\label{I:eq:alp-0}
			|b_{i,t}|
			=o\!\left((\log\lambda_t)^{-1/2}\right)
			=o\!\left(|\log\eps|^{-1/2}\right).
		\end{equation}
	\end{lemma}
	
	\begin{proof}
		Recall $\bar z=(\bar p,\bar\ell,0)$. Set
		$\Theta(p,\ell,b):=U_{p,\lambda(\ell),b}$ and use the coordinates
		\[
		\xi_{\bar z}(z):=
		\left(
		(\bar\lambda_t y_{\bar p_t}(p_t))_{t\in I'},\,
		(\ell_t-\bar\ell_t)_{t\in I'},\,
		(\sqrt{\bar L}\,b_{i,t})_{i\in I,t\in I'}
		\right).
		\]
		For sufficiently small $\eps_0$, the closure of
		$\mathcal D_{\bar z}(r_*)$ is compact, its centers lie in
		$M_{p^*,\tau}$, and they remain uniformly separated.
		
		Lemmas~\ref{I:lem:uniform-C2-bubble-map} and \ref{I:lem:tangent-gram}
		give
		\[
		\|D_\xi\Theta\,h\|_{\mathcal H}^2\ge c|h|^2,
		\qquad
		\|D_\xi\Theta\|+\|D_\xi^2\Theta\|\le C
		\]
		on a fixed small coordinate neighborhood containing the box.
		Indeed, the normalized amplitude, scale, and position directions are
		\[
		\bar L^{-1/2}X^b_{i,t},\qquad
		\lambda_tX_t^\lambda,\qquad
		\bar\lambda_t^{-1}X^p_{t,l},
		\quad l=1,2,
		\]
		up to the uniformly bounded coordinate changes in the position
		variables. Their limiting Gram matrix is positive definite.
		
		After decreasing $r_*$, Taylor's formula yields
		\begin{equation}\label{I:eq:modulation-bilipschitz}
			c|\xi_{\bar z}(z)-\xi_{\bar z}(z')|
			\le \|\Theta(z)-\Theta(z')\|_{\mathcal H}
			\le C|\xi_{\bar z}(z)-\xi_{\bar z}(z')|
		\end{equation}
		for $z,z'\in\overline{\mathcal D_{\bar z}(r_*)}$.
		The function $z\mapsto\|u-\Theta(z)\|_{\mathcal H}^2$ attains its
		minimum on the  compact set  $\overline{\mathcal D_{\bar z}(r_*)}$. Since
		\[
		\|u-\Theta(\bar z)\|_{\mathcal H}\le C\eps,
		\qquad
		\inf_{z\in\partial\mathcal D_{\bar z}(r_*)}
		\|u-\Theta(z)\|_{\mathcal H}\ge cr_*-C\eps,
		\]
		every minimizer is interior for sufficiently small $\eps_0$.
		Its Euler equations are
		\begin{align*}
			\langle u-\Theta,X^b_{i,t}\rangle_{\mathcal H}&=0,
			&&i\in I,\ t\in I',\\
			\langle u-\Theta,X_t^\lambda\rangle_{\mathcal H}&=0,
			&&t\in I',\\
			\langle u-\Theta,X^p_{t,l}\rangle_{\mathcal H}&=0,
			&&t\in I',\ l=1,2.
		\end{align*}
		
		For $F(\xi):=\frac12\|u-\Theta(\xi)\|_{\mathcal H}^2$,
		\[
		\begin{aligned}
			D^2F(\xi)[h,h]
			&=\|D_\xi\Theta\,h\|_{\mathcal H}^2
			-\langle u-\Theta,D_\xi^2\Theta[h,h]\rangle_{\mathcal H} \ge \bigl(c-C(r_*+\eps)\bigr)|h|^2.
		\end{aligned}
		\]
		Choosing $r_*$ and then $\eps_0$ small makes $F$ strictly convex
		on a convex coordinate neighborhood containing the box. Hence the
		interior minimizer is unique.
		
		By minimality and \eqref{I:eq:modulation-bilipschitz},
		\[
		\|u-\Theta(z)\|_{\mathcal H}\le C\eps,\qquad
		|\xi_{\bar z}(z)|
		\le C\|\Theta(z)-\Theta(\bar z)\|_{\mathcal H}
		\le C\eps.
		\]
		Now let $\bar z'$ be another preliminary ordered representation, with
		minimizer $z'$. Then
		\[
		\|\Theta(\bar z)-\Theta(\bar z')\|_{\mathcal H}\le C\eps.
		\]
		After relabeling, Lemma~\ref{I:lem:preliminary-identification} gives
		$
		|\xi_{\bar z}(\bar z')|\le C\eps.
		$
		Applying the preceding estimates in both charts and using their
		uniformly equivalent normalized coordinates, we obtain
		\[
		|\xi_{\bar z}(z)|+|\xi_{\bar z}(z')|\le C\eps.
		\]
		Thus $z$ and $z'$ lie in $\mathcal D_{\bar z}(r_*)$ and satisfy
		the same Euler equations. Strict convexity gives $z=z'$.
		
		Finally, $|\xi_{\bar z}(z)|\le C\eps$ implies
		\[
		\sqrt{\bar L}\max_{i,t}|b_{i,t}|\le C\eps.
		\]
		Since
		$
		\bar L\asymp\log\lambda_t,
		\log\lambda_t\ge|\log\eps|-r_*,
		$
		we conclude that
		\[
		\sqrt{\log\lambda_t}\,|b_{i,t}|\le C\eps,
		\qquad
		|b_{i,t}|=o\!\left(|\log\eps|^{-1/2}\right),
		\]
		which proves \eqref{I:eq:alp-0}.
	\end{proof}

	Explicitly, with
	$\alpha_{i,j,t}:=a_{ij}+\delta_{ij}b_{i,t}$, we have the following  decomposition: 
	\[
	u_i=\sum_{t\in I'}\sum_{j\in I}
	\alpha_{i,j,t}PW_{j,p_t,\lambda_t}+w_i,
	\qquad \forall i\in I,
	\]
with 
	\begin{align*}
		\begin{cases}
			\max_{i\in I}\|w_i\|< C\eps \quad\text{ and }\quad \max_{i\in I, t\in I'} |b_{i,t}| < \frac{C}{\sqrt{|\log \eps|}}, \\
			w:=( w_1,\ldots, w_n) \in E^{N}_{p,\lam,b},
		\end{cases}
	\end{align*}
	where the normal space is
	\[
	\begin{aligned}
		E^N_{p,\lambda,b}:=\Bigl\{w\in\mathcal H:\;&
		\langle w_i,PW_{i,p_t,\lambda_t}\rangle=0
		&&\forall i\in I,\ t\in I',\\
		&\sum_{i,j\in I}\alpha_{i,j,t}
		\langle w_i,\partial_{\lambda_t}PW_{j,p_t,\lambda_t}\rangle=0
		&& \forall t\in I',\\
		&\sum_{i,j\in I}\alpha_{i,j,t}
		\langle w_i,\partial_{(p_t)_l}PW_{j,p_t,\lambda_t}\rangle=0
		&& \forall t\in I',\ l=1,2
		\Bigr\}.
	\end{aligned}
	\]

	We define these moments for amplitudes satisfying
	\[
	|b_{i,t}|\le b_0,
	\qquad
	0<b_0<\min_{i\in I}\frac{q_i}{2\sigma_i}.
	\]
	Then $q_i+\sigma_i b_{i,t}\ge q_i/2>0$ for every
	$i\in I$ and $t\in I'$. Since
	\[
	\sum_j\alpha^{i,j,t}v_j(y)
	=-(m_i+\sigma_i b_{i,t})\log|y|+O(1),
	\]
	the moments below are finite. 
	
	Recall that $\sigma_i=(2\pi)^{-1}\int_{\R^2}e^{v_i}$. For $i,j\in I$, $t\in I'$, and $\alpha_{i,j,t}=a_{ij}+\delta_{ij}b_{i,t}$, we introduce only the additional profile moments needed below:
	\begin{align*}
		\sigma_{i,t}^{\alpha}:= \frac 1 {2\pi} \int_{\R^2} e^{ \sum_{j\in I} \alpha^{i,j,t} v_j(y)}\, dy;\\
		\sigma^{(1)}_{i,j}:=  \frac 1 {2\pi}\int_{\R^2} e^{v_i(y)} v_j(y)\, dy;\qquad &
		\sigma^{(1,\alpha)}_{i,j,t}:=  \frac 1 {2\pi}\int_{\R^2} e^{\sum_{\ell\in I}\alpha^{i,\ell,t}v_\ell(y)}v_j(y)\,dy;  \\
		\sigma_{i,j}^{(2)}:= \frac 1 {2\pi} \int_{\R^2}  e^{v_i(y)} \nabla v_j  (y)\cdot y \, dy; \qquad &
		\sigma^{(2,\alpha)}_{i,j, t}:=\frac 1 {2\pi }\int_{\R^2}   e^{\sum_{\ell \in I}\alpha^{i, \ell , t}  v_\ell (y)}  \nabla v_{j}(y)\cdot y \, dy; \\
		\sigma^{\ln }_{i}:=  \frac 1 {2\pi}\int_{\R^2} e^{v_i(y)} \ln |y|\, dy;\qquad &
		\sigma^{\ln,\alpha }_{i,t}:=  \frac 1 {2\pi}\int_{\R^2} e^{\sum_{\ell\in I}\alpha^{i,\ell,t}v_\ell(y)}\ln|y|\,dy.
	\end{align*}

	%-------------------------Section 4
	\subsection{Normalized equations and local masses}
	\label{I:sec:coarse-blowup}

	Let $\rho^k\to\rho^*\in\Gamma_N$, where either $\rho^k\equiv\rho^*$
	or $\{\rho^k\}$ is admissible as in
	Definition~\ref{I:def:admissible-path}. Let $u^k\in\mathcal H$ solve
	\begin{equation}\label{I:eq:solution}
		-\Delta_g u_i^k
		=\sum_{j=1}^n a_{ij}\rho_j^k
		\left(
		\frac{h_je^{u_j^k}}{\int_M h_je^{u_j^k}\,dV_g}-1
		\right),
		\qquad i\in I.
	\end{equation}
	Following \cite[Section~2]{HuangZhang2022}, define
	\begin{equation*}
		%\label{I:eq:normalized-solution}
		\Theta_i^k:=
		u_i^k-\log\int_M h_ie^{u_i^k}\,dV_g.
	\end{equation*}
	Then
	\begin{equation}\label{I:eq:normalized-equation}
		-\Delta_g\Theta_i^k
		=\sum_{j=1}^n a_{ij}\rho_j^k
		(h_je^{\Theta_j^k}-1),
		\qquad
		\int_M h_ie^{\Theta_i^k}\,dV_g=1.
	\end{equation}
	
	Assume
	\[
	\max_{i\in I}\|u_i^k\|_{L^\infty(M)}\to\infty.
	\]
	Under \ref{I:assump:H1}--\ref{I:assump:H2}, the blow-up analysis
	of Lin and Zhang \cite{LinZhang2011} gives, after a subsequence,
	exactly $N$ distinct concentration points at which all components
	blow up:
	\[
	p^*:=(p_1^*,\ldots,p_N^*)\in M^N\setminus\Delta.
	\]
	Choose mutually disjoint geodesic disks $B_{2r_0}(p_t^*)$.
	Let $\widehat p_t^k\in B_{r_0}(p_t^*)$ be the refined local maximum
	points supplied by the local approximation, and set
	\begin{equation}\label{I:eq:height-scale}
		M_{k,t}:=\max_{i\in I}
		\left\{
		\Theta_i^k(\widehat p_t^k)
		+\log\bigl(\rho_i^kh_i(\widehat p_t^k)\bigr)
		\right\}.
	\end{equation}
	Define
	\begin{equation*}
		%\label{I:eq:epsilon-lambda-relation}
		\varepsilon_{k,t}:=e^{-M_{k,t}/2},
		\qquad
		\widehat\lambda_t^k:=\varepsilon_{k,t}^{-1}.
	\end{equation*}
	Thus
	$
	M_{k,t}-M_{k,s}
	=2\log\frac{\widehat\lambda_t^k}{\widehat\lambda_s^k},
	$
	so bounded height differences are equivalent to comparable scales.
	Set 
	$
	\Psi_k:=\Psi_{\widehat m}(\widehat\lambda_1^k).
	$
	
	Let $\{\Omega_t^k\}_{t=1}^N$ be a measurable partition of $M$,
	up to boundaries of zero area, with each $\Omega_t^k$ containing
	a fixed neighborhood of $p_t^*$. Define
	\[
	\mathfrak p_{i,t}^k
	:=\rho_i^k\int_{\Omega_t^k}h_ie^{\Theta_i^k}\,dV_g,
	\qquad
	\mu_{i,t}^k:=\frac{\mathfrak p_{i,t}^k}{2\pi}.
	\]
	Then
	\[  \sum_{t\in I'} \mathfrak p_{i,t}^k=\rho_i^k,\qquad
	\sum_{t\in I'}\mu_{i,t}^k=\frac{\rho_i^k}{2\pi}. \]

	%--------------
	
	\subsection{Refined local blow-up estimates}
	
	We record the local estimates of
	\cite[Theorems~4.2--4.3]{LinZhang2013} used to expand the masses
	and moments in Section~\ref{I:sec:integrated-critical-estimates}.
	For $\widehat m=4$, the additional mass hypothesis follows from
	Theorem~\ref{I:thm:path-free-coarse-matching}, proved using the
	first-order local estimate and the mass variance identity.
	
	\begin{proposition}
		\label{I:prop:refined-local-Harnack}
		Assume (H1)--(H2), (A1), and $h_i\in C^3(M)$ with $h_i>0$
		for every $i\in I$.
		Let $u^k$ be a sequence of blow-up solutions with
		$\rho^k\to\rho^*\in\Gamma_N$. Assume either
		$\rho^k\equiv\rho^*$ or $\{\rho^k\}$ is admissible in
		Definition~\ref{I:def:admissible-path}, with
		$\Lambda_{I,N}(\rho^k)>0$ when $N>1$.  For $N=1$, either transverse side is allowed.
		
		Fix a concentration point $p_t^*$. In the chart centered at
		$\widehat p_t^k$, choose $f_{i,t}^k$ satisfying
		\[
		\Delta f_{i,t}^k
		=(A\rho^k)_i e^{\widehat\varphi_{\widehat p_t^k}},
		\qquad
		f_{i,t}^k(0)=0,
		\qquad
		\|f_{i,t}^k\|_{C^3}\le C.
		\]
		Let $\phi_{i,t}^k$ be harmonic, with $\phi_{i,t}^k(0)=0$,
		such that $\Theta_i^k-f_{i,t}^k-\phi_{i,t}^k$ is constant
		on the boundary of the coordinate disk. For
		$y=y_{\widehat p_t^k}(x)/\varepsilon_{k,t}$, set
		\[
		v_{i,t}^k(y):=
		\Theta_i^k(x)
		+\log\bigl(\rho_i^kh_i(\widehat p_t^k)\bigr)
		+2\log\varepsilon_{k,t}
		-f_{i,t}^k(x)-\phi_{i,t}^k(x).
		\]
		
		There exist a radial entire solution $V^{k,t}$ of
		the PDE in \eqref{I:eq:entire}, with its own mass vector, a displacement
		$\zeta_t^k=O(\varepsilon_{k,t})$, and an angular correction
		\[
		\mathcal A_{i,t}^k(r,\theta)
		:=\varepsilon_{k,t}
		\bigl(G_{1,i,t}^k(r)\cos\theta
		+G_{2,i,t}^k(r)\sin\theta\bigr)
		\]
		with the following properties. Denote
		\[
		m_{i,t}^{V,k}
		:=\frac1{2\pi}\sum_j a_{ij}
		\int_{\mathbb R^2}e^{V_j^{k,t}}\,dy,
		\qquad
		\widehat m_{k,t}:=\min_i m_{i,t}^{V,k}\to\widehat m.
		\]
		For $|\alpha|\le1$ and $|y|\le c\varepsilon_{k,t}^{-1}$,
		\[
		\begin{aligned}
			&\quad \left|
			D^\alpha\left[
			v_{i,t}^k
			-V_i^{k,t}(\,\cdot-\zeta_t^k)
			-\mathcal A_{i,t}^k(\,\cdot-\zeta_t^k)
			\right](y)
			\right|\\
			&\qquad\le C\varepsilon_{k,t}^2
			\begin{cases}
				(1+|y|)^{4-\widehat m_{k,t}-|\alpha|}
				\log(2+|y|),&\widehat m<4,\\
				(1+|y|)^{-|\alpha|}
				\bigl(\log(2+|y|)\bigr)^2,&\widehat m=4.
			\end{cases}
		\end{aligned}
		\]
		Moreover, for $l=1,2$ and $0\le r\le c\varepsilon_{k,t}^{-1}$,
		\[
		|G_{l,i,t}^k(r)|
		\le Cr
		\begin{cases}
			(1+r)^{2-\widehat m_{k,t}},&\widehat m<4,\\
			(1+r)^{-2},&\widehat m=4.
		\end{cases}
		\]
		When $\widehat m=4$, we also have
		\[
		\max_i|m_{i,t}^{V,k}-4|
		\le \frac{C}{\log\varepsilon_{k,t}^{-1}}.
		\]
		All constants are independent of $k,i,t$.
	\end{proposition}
	\begin{proof}
		Lemma~\ref{I:lem:path-free-local-conversion} verifies the local
		hypotheses (3.1)--(3.6) of \cite{LinZhang2013}. The harmonic
		function $\phi_{i,t}^k$ is defined as in (3.9) of that paper.
		
		By \cite[Lemma~4.1]{LinZhang2013}, the maximum after harmonic
		subtraction has scaled displacement
		$\zeta_t^k=O(\varepsilon_{k,t})$. The second-order expansion
		uses the same radial solution $V^{k,t}$, translated by
		$\zeta_t^k$, together with the first angular correction.
		In particular, its masses are unchanged:
		\[
		\int_{\mathbb R^2}e^{V_i^{k,t}(y-\zeta_t^k)}\,dy
		=\int_{\mathbb R^2}e^{V_i^{k,t}(y)}\,dy.
		\]
		For $\widehat m<4$, the estimates therefore follow from
		\cite[Theorem~4.2, (4.49)--(4.50)]{LinZhang2013}.
		
		For $\widehat m=4$, set $L_k:=\frac1N\sum_t\log\varepsilon_{k,t}^{-1}.
		$
		Theorem~\ref{I:thm:path-free-coarse-matching} gives, for some
		fixed $\beta>0$,
		\[
		\max_{i,t}|m_{i,t}^{V,k}-4|\le Ce^{-\beta L_k},
		\qquad
		\log\varepsilon_{k,t}^{-1}=L_k+O(1).
		\]
		Hence
		\[
		\max_{i,t}
		\left\{|m_{i,t}^{V,k}-4|
		\log\varepsilon_{k,t}^{-1}\right\}
		\le CL_ke^{-\beta L_k}\longrightarrow0.
		\]
		This verifies the additional hypothesis of
		\cite[Theorem~4.3]{LinZhang2013}; its (4.55)--(4.56)
		give the remaining estimates. A fixed rescaling of the coordinate
		disk changes only the constants.
	\end{proof}
	
	%--------------------------------

	\section{Mass matching and sharp parameter estimates}\label{I:sec:sharp-defect}
	
	\subsection{Mass variance and scale balance}
	
	Set
	\[
	P(a):=4\mathbf1^Ta-a^TAa,
	\qquad q:=A\sigma-2\mathbf1,
	\qquad \sigma=\frac{\rho^*}{2\pi N}.
	\]
	Here,  $\mathbf1:=(1,\ldots,1)^T\in\mathbb R^n$.
	Then $P(\sigma)=0$ and $q_i=m_i-2>0$.
	
	\begin{lemma}
		Under (H1)--(H2), $A$ and $A^{-1}$ have one positive and $n-1$ negative
		eigenvalues. Moreover,
		\begin{equation}\label{I:eq:radial-timelike}
			q^TA^{-1}q=4\mathbf1^TA^{-1}\mathbf1>0,
		\end{equation}
		and, for some $c>0$,
		\begin{equation}\label{I:eq:radial-tangent-coercivity}
			x^TAx\le-c|x|^2\qquad(q^Tx=0).
		\end{equation}
	\end{lemma}
	
	\begin{proof}
		The symmetric matrix $A^{-1}$ is irreducible and has nonnegative off-diagonal
		entries. Perron--Frobenius gives a simple largest eigenvalue $\nu$ with
		eigenvector $v>0$. Since $A^{-1}\mathbf1\ge0$ is nonzero,
		\[
		\nu\ge\frac{\mathbf1^TA^{-1}\mathbf1}{n}>0,
		\qquad Av=\nu^{-1}v.
		\]
		Thus $\nu^{-1}$ is the spectral radius of $A$. Any other positive
		eigenvalue $\zeta<\nu$ of $A^{-1}$ would give an eigenvalue
		$\zeta^{-1}>\nu^{-1}$ of $A$, a contradiction. This proves the eigenvalue assertion.
		Expanding $q=A\sigma-2\mathbf1$ and using $P(\sigma)=0$ proves
		\eqref{I:eq:radial-timelike}. The vector $A^{-1}q$ satisfies
		\[
		(A^{-1}q)^TA(A^{-1}q)>0,\qquad (A^{-1}q)^TAx=q^Tx.
		\]
		Since $A$ has exactly one positive eigenvalue and
		$(A^{-1}q)^TA(A^{-1}q)>0$, the quadratic form $x^TAx$ is negative definite
		on
		\[
		\{x:(A^{-1}q)^TAx=0\}=q^\perp.
		\]
		Compactness of the unit sphere in $q^\perp$ gives
		\eqref{I:eq:radial-tangent-coercivity}.
	\end{proof}
	
	\begin{proposition}\label{I:prop:radial-variance}
		Let $P(a_t)=0$ for $t=1,\ldots,N$, and set
		\[
		\bar a:=\frac1N\sum_ta_t,\qquad
		d_t:=a_t-\bar a,\qquad
		e:=a-\bar a.
		\]
		Then
		\begin{equation}\label{I:eq:radial-variance-exact}
			P(a)=\frac1N\sum_td_t^TAd_t
			+(4\mathbf1-2A\bar a)^Te-e^TAe.
		\end{equation}
		If $\max_t|a_t-\sigma|$ is sufficiently small, then
		\begin{equation}\label{I:eq:radial-variance-negative}
			\frac1N\sum_td_t^TAd_t\le-\frac cN\sum_t|d_t|^2.
		\end{equation}
		Consequently, $P(a)\ge 0$ implies 
		\begin{equation}\label{I:eq:radial-variance-bound}
			\frac1N\sum_t|d_t|^2\le C(|e|+|e|^2).
		\end{equation}
	\end{proposition}
	
	\begin{proof}
		Averaging $P(\bar a+d_t)=0$ gives
		$P(\bar a)=N^{-1}\sum_td_t^TAd_t$; expansion at $\bar a$ proves
		\eqref{I:eq:radial-variance-exact}.
		Write $u_t=a_t-\sigma$ and $\max_t|u_t|\le\varepsilon$. Subtracting
		$P(\sigma+u_t)=P(\sigma+u_s)=0$ gives
		\[
		|q^T(u_t-u_s)|\le C\varepsilon|u_t-u_s|.
		\]
		For $n_t=(q^Td_t)q/|q|^2$ and $w_t=d_t-n_t$, Jensen's inequality gives
		\[
		\sum_t|n_t|^2\le C\varepsilon^2\sum_t|d_t|^2,
		\qquad
		d_t^TAd_t\le-\frac c2|w_t|^2+C|n_t|^2.
		\]
		Summation and smallness of $\varepsilon$ prove
		\eqref{I:eq:radial-variance-negative}. Inserting this bound into
		\eqref{I:eq:radial-variance-exact} proves
		\eqref{I:eq:radial-variance-bound}.
	\end{proof}
	
	\begin{lemma}
		For $z\in\mathbb R^n$ and $r\in\mathbb R$,
		\begin{equation}\label{I:eq:full-linear-matching-inverse}
			|z|+|r|\le C\bigl(|Az+qr|+|q^Tz|\bigr).
		\end{equation}
	\end{lemma}
	
	\begin{proof}
		If $Az+qr=0$ and $q^Tz=0$, then
		$r q^TA^{-1}q=0$, hence $r=0$ and $z=0$ by
		\eqref{I:eq:radial-timelike}. The linear map
		$(z,r)\mapsto(Az+qr,q^Tz)$ therefore has a bounded inverse.
	\end{proof}
	
	For a separated configuration $p$, write
	\[
	H_{i,t}(p):=\log\rf_{p,t}^i(p_t),\qquad
	H_t(p):=(H_{i,t}(p))_{i\in I},\qquad
	\overline H(p):=\frac1N\sum_tH_t(p).
	\]
	
	\begin{lemma}\label{I:lem:stationary-profile-matching}
		Assume (H1)--(H2) and (A1). Let $L_k\to\infty$,
		$p^k\to p^*\in M^N\setminus\Delta$, and let $\delta^k\in E_0$ be bounded.
		Suppose $P(a_t^k)=0$ and $a_t^k\to\sigma$. Set
		\[
		\bar a^k:=\frac1N\sum_ta_t^k,\qquad
		d_t^k:=a_t^k-\bar a^k,\qquad
		D_k^2:=\frac1N\sum_t|d_t^k|^2.
		\]
		Let $E_k>0$ satisfy $L_k^2E_k\to0$. If
		\begin{align}
			|\bar a^k-\sigma|&\le CE_k,
			\notag\\
			\max_t|L_kAd_t^k+q\delta_t^k+H_t(p^k)-\overline H(p^k)|
			&\le C(D_k+L_kE_k),
			\label{I:eq:stationary-explicit-matching-system}\\
			|\nabla\ff(p^k)|&\le C(D_k+E_k),
			\label{I:eq:stationary-local-location-input}
		\end{align}
		then, after a subsequence, $(p^*,\delta^*)$ is a balanced critical pair and
		\[
		|p^k-p^*|+\max_t|a_t^k-\sigma|\le CE_k,
		\qquad |\delta^k-\delta^*|\le CL_kE_k.
		\]
		For the corresponding maximum-zero radial profiles and each fixed $R$,
		\begin{equation}\label{I:eq:stationary-matching-lemma-profile}
			\max_{s,t}\|V^{k,t}-V^{k,s}\|_{C^2(B_R)}\le C_RE_k.
		\end{equation}
	\end{lemma}
	
	\begin{proof}
		Apply \eqref{I:eq:radial-variance-bound} with
		$a=\sigma$ and $e=\sigma-\bar a^k$. Then
		\[
		D_k^2\le CE_k,\qquad L_kD_k\to0.
		\]
		Passing to a subsequence in \eqref{I:eq:stationary-explicit-matching-system}
		and \eqref{I:eq:stationary-local-location-input} gives
		\[
		\nabla\ff(p^*)=0,\qquad
		H_t(p^*)-\overline H(p^*)+q\delta_t^*=0.
		\]
		By (A1), $|p^k-p^*|\le C(D_k+E_k)$. Subtract the limiting balance
		equation and use the relative modes $e^\alpha$:
		\[
		A(L_kd^{k,\alpha})+q(r_\alpha^k-r_\alpha^*)
		=O(D_k+L_kE_k),\qquad
		d^{k,\alpha}:=\sum_te_t^\alpha d_t^k.
		\]
		Expanding $P(a_t^k)=0$ at $\bar a^k$ also gives
		\[
		(A\bar a^k-2\mathbf1)^Td_t^k
		=\frac12\left(\frac1N\sum_s(d_s^k)^TAd_s^k
		-(d_t^k)^TAd_t^k\right),
		\]
		so $|q^Td_t^k|\le C(D_k^2+E_kD_k)$. Therefore
		\eqref{I:eq:full-linear-matching-inverse} yields
		\[
		L_kD_k+|r^k-r^*|
		\le C(D_k+L_kE_k+L_kD_k^2+L_kE_kD_k).
		\]
		Absorbing the terms containing $D_k$ gives
		$D_k\le CE_k$ and $|r^k-r^*|\le CL_kE_k$.
		The location and absolute mass bounds follow.
		Finally, Lemma~\ref{I:lem:radial-family-derivatives} and the locally
		Lipschitz maximum-zero normalization give
		\eqref{I:eq:stationary-matching-lemma-profile}.
	\end{proof}
	
	\subsection{Local estimates and comparable concentration scales}
	
	Write $s_i^k:=\rho_i^kh_ie^{\Theta_i^k}$. By
	\cite[Proposition~4.1 and the proof of Theorem~1.1]{LinZhang2011},
	the concentration points are distinct, their limiting mass vectors equal
	$\sigma$, the rescaled limits are full, and the exterior mass tends to zero.
	Componentwise oscillations on fixed circles are bounded. These conclusions
	hold for both stationary and moving sequences.
	
	Fix disjoint chart disks of radius $8r_0$, with the matching circles
	contained where the cutoffs equal $1$. Constants depend on the fixed
	data, $r_0$, the separated center set, the compact mass neighborhood,
	the positive lower bounds and $C^3$ norms of $h_i$, and, when used,
	the Morse inverse bound and transverse direction.
	Set
	\[
	q_*:=\widehat m-2,\qquad
	0<\beta<\tfrac18\min\{1,q_*\}.
	\]
	
	\begin{lemma}
		\label{I:lem:path-free-local-conversion}
		There are centers $p_t^k\to p_t^*$, lengths $\epsilon_t^k\downarrow0$,
		and radial entire solutions $V^{k,t}$ of
		$-\Delta V_i^{k,t}=\sum_j a_{ij}e^{V_j^{k,t}}$, centered at the origin and normalized by  $\max_{i,t} V^{k,t}_i(0)=0$, such that, on $|x|<2r_0$ in the centered chart,
		\begin{equation}\label{I:eq:path-free-local-form}
			\begin{aligned}
				\Theta_i^k(x)+\log(\rho_i^kh_i(p_t^k))-f_{i,t}^k(x)
				&=2\ell_t^k+V_i^{k,t}(x/\epsilon_t^k)\\
				&\quad
				+\phi_{i,t}^k(x)+\xi_{i,t}^k(x),
				\qquad \ell_t^k=\log(\epsilon_t^k)^{-1}.
			\end{aligned}
		\end{equation}
		Here $\Delta f_{i,t}^k=(A\rho^k)_i e^{\widehat\varphi_{p_t^k}}$,
		$f_{i,t}^k(0)=0$, $\phi_{i,t}^k$ is harmonic with
		$\phi_{i,t}^k(0)=0$, and $f_{i,t}^k,\phi_{i,t}^k$ have uniformly bounded
		$C^2$ norms on this disk. 
		For every fixed $0<r<R<2r_0$, define
		\[
		\mathcal A_{r,R}:=\{x\in\mathbb R^2:r\le |x|\le R\}.
		\]
		Then
		\begin{equation}\label{I:eq:path-free-first-order-bound}
			\|\xi_{i,t}^k\|_{L^\infty(B_{2r_0})}
			+\|\nabla_x\xi_{i,t}^k\|_{L^\infty(\mathcal A_{r,R})}
			\le C_{r,R}(\epsilon_t^k)^{2\beta},
		\end{equation}
		where $C_{r,R}$ is independent of $k,i,t$.
		Set
		$
		a_{i,t}^k=\frac1{2\pi}\int_{\R^2}e^{V_i^{k,t}},\, 
		m_{i,t}^{V,k}=(Aa_t^k)_i.
		$
		Then $P(a_t^k)=0$, $a_t^k\to\sigma$, and, for fixed $r\in(0,r_0]$,
		\begin{align}
			\left|\frac1{2\pi}\int_{U_r(p_t^k)}s_i^k\,dV_g-a_{i,t}^k\right|
			&\le C_r(\epsilon_t^k)^{2\beta},
			\label{I:eq:path-free-disk-mass}\\
			\int_{U_r(p_t^k)}d_g(x,p_t^k)s_i^k\,dV_g
			&\le C_r(\epsilon_t^k)^{2\beta}.
			\label{I:eq:path-free-first-moment}
		\end{align}
	\end{lemma}
	
	\begin{proof}
		We apply \cite[Theorem~3.1, (3.11)]{LinZhang2013}.
		Solve $\Delta f_{i,t}^k=(A\rho^k)_i e^{\widehat\varphi_{p_t^k}}$
		on a larger fixed chart disk with zero boundary values, and subtract
		its value at the center. Interior elliptic estimates give uniform
		$C^3$ bounds on the smaller disk. Set
		\[
		\begin{aligned}
			w_i^k&:=\Theta_i^k+\log(\rho_i^kh_i(p_t^k))-f_{i,t}^k,\qquad 
			K_i^k(x):=\frac{h_i(x)}{h_i(p_t^k)}
			e^{\widehat\varphi_{p_t^k}(x)+f_{i,t}^k(x)}.
		\end{aligned}
		\]
		Equation~\eqref{I:eq:normalized-equation} gives
		\[
		-\Delta w_i^k=\sum_j a_{ij}K_j^ke^{w_j^k},\qquad
		K_i^ke^{w_i^k}\,dx=s_i^k\,dV_g,
		\]
		and
		\[
		K_i^k(0)=1,\qquad
		C^{-1}\le K_i^k\le C,\qquad
		\|K_i^k\|_{C^3}\le C,\qquad
		\int_{B_{2r_0}}e^{w_i^k}\,dx\le C\rho_i^k\le C.
		\]
		The blow-up alternative gives one isolated concentration point,
		bounded componentwise boundary oscillations, and a full rescaled limit.
		Subtracting $f_{i,t}^k$ preserves these properties. After a fixed disk
		rescaling, these are precisely (3.1)--(3.6) of the cited theorem.
		
		For the centering, first choose a maximum of one component of the
		background-subtracted solution on a smaller fixed disk, and relabel that
		component as $1$. Recenter and normalize the conformal chart there; a
		constant can be subtracted from each $f_i^k$ to impose $f_i^k(0)=0$.
		Subtract the harmonic
		extension of the boundary oscillation; it has value zero at the center
		and uniformly bounded derivatives on the smaller disk. In the source
		theorem the radial solution agrees with all component values at the
		center, and component $1$ fixes the dilation. Suppress $k,t$ and write
		$\widetilde V$, $\widetilde\epsilon$ for those source parameters. Set
		\[
		c:=\max_i\widetilde V_i(0),\qquad
		V_i(y):=\widetilde V_i(e^{-c/2}y)-c,\qquad
		\epsilon:=e^{-c/2}\widetilde\epsilon.
		\]
		Then the entire masses are unchanged and
		\[
		\max_i V_i(0)=0,\qquad
		2\log\epsilon^{-1}+V_i(x/\epsilon)
		=2\log\widetilde\epsilon^{-1}
		+\widetilde V_i(x/\widetilde\epsilon).
		\]
		Since $f_i^k(0)=\phi_i^k(0)=0$ and the radial values match at the center,
		$2\log\epsilon^{-1}$ is the height in
		\eqref{I:eq:height-scale} at that center. The maximum-zero profile
		family is locally Lipschitz in its masses by
		Lemma~\ref{I:lem:radial-family-derivatives} and this dilation.
		On a fixed rescaled neighborhood, the radial Hessian is uniformly
		negative definite, whereas subtraction of a smooth background changes
		the scaled gradient by $O(\epsilon_t^k)$. Hence the physical center
		and height change by $O((\epsilon_t^k)^2)$, which is absorbed in
		\eqref{I:eq:path-free-first-order-bound}.
		
		Estimate (3.11), for small fixed $\eta>0$ and $|\alpha|\le1$, is
		\[
		|D_y^\alpha\xi_i(\epsilon y)|\le
		\begin{cases}
			C\epsilon(1+|y|)^{3-\widehat m+\eta-|\alpha|},&\widehat m\le3,\\
			C\epsilon(1+|y|)^{\eta-|\alpha|},&\widehat m>3.
		\end{cases}
		\]
		Choose $2\beta<\min\{q_*-\eta,1-\eta\}$.
		Taking the supremum and converting the first derivative on a fixed
		annulus prove \eqref{I:eq:path-free-first-order-bound}.
		
		The radial ODE, uniformly on a small compact mass neighborhood, gives
		\[
		e^{V_i^{k,t}(y)}\le C(1+|y|)^{-2-q_*/2},\quad
		V_i^{k,t}(y)=-m_{i,t}^{V,k}\log|y|+c_i(a_t^k)
		+O(|y|^{-q_*/2}).
		\]
		Lemma~\ref{I:lem:radial-family-derivatives} gives the local Lipschitz
		dependence of $c_i$. The entire Pohozaev identity gives $P(a_t^k)=0$;
		local convergence and the uniform tail bound give $a_t^k\to\sigma$.
		
		On changing variables in the source measure, put
		$J_i^k(x)=K_i^k(x)e^{\phi_i^k(x)}$. Then $J_i^k(0)=1$,
		$|J_i^k(x)-1|\le C|x|$, and the rescaled density is
		$J_i^k(\epsilon y)e^{V_i(y)+\xi_i(\epsilon y)}$. The elementary bounds
		\[
		\int_{|y|>r/\epsilon}e^{V_i(y)}dy\le C_r\epsilon^{q_*/2},\qquad
		\epsilon\int_{|y|<r/\epsilon}|y|e^{V_i(y)}dy
		\le C_r\epsilon^{2\beta}
		\]
		follow by integrating $(1+|y|)^{-2-q_*/2}$; the strict inequality
		$2\beta<\min\{1,q_*/2\}$ also absorbs a possible logarithm.
		Together with the supremum bound for $\xi$, these prove both assertions.
	\end{proof}
	
	\begin{lemma}
		\label{I:lem:path-free-circle}
		Let $\ell_-^k=\min_t\ell_t^k$, $\ell_+^k=\max_t\ell_t^k$, and
		$E_0^k=e^{-2\beta\ell_-^k}$. Let $\Omega_t^k$ be a measurable partition
		containing the disjoint disks $U_{r_0}(p_t^k)$, with boundaries of zero
		area, and set $\mu_{i,t}^k=(2\pi)^{-1}\int_{\Omega_t^k}s_i^kdV_g$.
		Then
		\begin{equation}\label{I:eq:path-free-mass-ledger}
			\max_{i,t}|\mu_{i,t}^k-a_{i,t}^k|\le CE_0^k,\quad
			\left|\bar a^k-\frac{\rho^k}{2\pi N}\right|\le CE_0^k,
			\quad \ell_+^k\le C(1+\ell_-^k).
		\end{equation}
		Define the regular potential at the $t$-th disk by
		\[
		S_{i,t}^k(x)=2\pi m_{i,t}^{V,k}H^g(x,p_t^k)
		+2\pi\sum_{s\ne t}m_{i,s}^{V,k}G^g(x,p_s^k).
		\]
		With $\overline\Theta_i^k=\int_M\Theta_i^kdV_g$, the circle identity is
		\begin{equation}\label{I:eq:path-free-circle-identity}
			\overline\Theta_i^k
			=-(m_{i,t}^{V,k}-2)\ell_t^k+c_i(a_t^k)
			-\log\rho_i^k-\log h_i(p_t^k)-S_{i,t}^k(p_t^k)+O(E_0^k).
		\end{equation}
		The error is uniform over the finitely many fixed matching circles.
	\end{lemma}
	
	\begin{proof}
		On a fixed annulus Lemma~\ref{I:lem:path-free-local-conversion} gives
		$\Theta_i^k=-(m_{i,t}^{V,k}-2)\ell_t^k+O(1)$. Since
		$m_{i,t}^{V,k}-2\ge q_*/2$ for large $k$, this is at most
		$-(q_*/2)\ell_-^k+C$. On this fixed exterior, the blow-up alternative gives
		$h_j e^{\Theta_j^k}=o(1)$ uniformly. Hence, for large $k$,
		\[
		\Delta_g\Theta_i^k
		=\sum_j a_{ij}\rho_j^k(1-h_j e^{\Theta_j^k})\ge0.
		\]
		The maximum principle therefore gives
		\begin{equation}\label{I:eq:path-free-exterior}
			\|s_i^k\|_{L^\infty(M\setminus\bigcup_tU_{r_0/2}(p_t^k))}
			\le Ce^{-(q_*/2)\ell_-^k}\le CE_0^k.
		\end{equation}
		For fixed $r\le r_0$, the exact mass decomposition is
		\begin{align*}
			2\pi(\mu_{i,t}^k-a_{i,t}^k)
			={}&\int_{B_{r/\epsilon_t^k}}
			\bigl(J_{i,t}^k(\epsilon_t^ky)e^{V_i^{k,t}(y)+\xi_{i,t}^k(\epsilon_t^ky)}
			-e^{V_i^{k,t}(y)}\bigr)dy\\
			&-\int_{\R^2\setminus B_{r/\epsilon_t^k}}e^{V_i^{k,t}(y)}dy
			+\int_{\Omega_t^k\setminus U_r(p_t^k)}s_i^k\,dV_g.
		\end{align*}
		The preceding lemma and \eqref{I:eq:path-free-exterior} bound every
		term by $CE_0^k$. Summing uses the exact equality
		$\sum_t\mu_t^k=\rho^k/(2\pi)$ and proves the first two claims of
		\eqref{I:eq:path-free-mass-ledger}.
		
		On every fixed matching annulus Green representation now gives
		\begin{equation}\label{I:eq:path-free-green-C1}
			\Theta_i^k(x)=\overline\Theta_i^k
			+2\pi\sum_sm_{i,s}^{V,k}G^g(x,p_s^k)+O_{C^1}(E_0^k).
		\end{equation}
		Here is a justification including the derivative. Split each source into
		a core disk of half the matching radius and its complement. In the
		core, Taylor's formula for $G^g(x,\cdot)$ and its first $x$ derivative
		has bounded coefficients, and its error is bounded by the first moment
		\eqref{I:eq:path-free-first-moment}. Replacing the core mass by
		$2\pi a_s^k$ uses \eqref{I:eq:path-free-disk-mass}. On the complement
		the density is $O(E_0^k)$; both the logarithmic kernel and its
		$d_g(x,y)^{-1}$ first derivative have bounded integrals on the surface.
		This proves \eqref{I:eq:path-free-green-C1}.
		
		Average this identity on a centered circle $|x|=r$. Its singular part is
		$-m_{i,t}^{V,k}\log r$. By \eqref{I:eq:path-free-local-form}, the other
		expression for the same average is
		\[
		-(m_{i,t}^{V,k}-2)\ell_t^k-m_{i,t}^{V,k}\log r+c_i(a_t^k)
		-\log(\rho_i^kh_i(p_t^k))+\frac1{2\pi r}\int_{\partial B_r}f_{i,t}^k+O(E_0^k).
		\]
		The harmonic $\phi_{i,t}^k$ has average zero. Furthermore
		\[
		\Delta(S_{i,t}^k-f_{i,t}^k)
		=\bigl(2\pi\sum_sm_{i,s}^{V,k}-(A\rho^k)_i\bigr)
		e^{\widehat\varphi_{p_t^k}}=O(E_0^k).
		\]
		Solving the zero-boundary Poisson equation on the fixed disk and using
		the mean-value property of its harmonic remainder shows that
		$\frac1{2\pi r}\int_{\partial B_r}(S_{i,t}^k-f_{i,t}^k)
		=S_{i,t}^k(0)+O(E_0^k)$. This proves
		\eqref{I:eq:path-free-circle-identity}, including its constant and sign.
		Subtract it for two indices $s,t$. All regular terms and $c_i(a_t^k)$
		are bounded, so
		$(m_{i,t}^{V,k}-2)\ell_t^k-(m_{i,s}^{V,k}-2)\ell_s^k=O(1)$.
		The positive lower and upper bounds for these exponents give
		$\ell_+^k\le C(1+\ell_-^k)$, without assuming comparable scales.
	\end{proof}
	
	\begin{lemma}
		\label{I:lem:path-free-location}
		In the notation above,
		\begin{equation}\label{I:eq:path-free-weighted-location}
			\sum_i a_{i,t}^k
			\left[\nabla\log h_i(p_t^k)+\nabla_xS_{i,t}^k(p_t^k)\right]
			=O(E_0^k).
		\end{equation}
		Consequently, writing $D_k^2=N^{-1}\sum_t|a_t^k-\bar a^k|^2$,
		\begin{equation}\label{I:eq:path-free-location-estimate}
			|\nabla\ff(p^k)|\le C\left(D_k+E_0^k+|\rho^k-\rho^*|\right).
		\end{equation}
	\end{lemma}
	
	\begin{proof}
		Fix a coordinate direction $l$ and a disk $B_r$ in the $t$-th chart.
		For $w_i=\Theta_i^k+\log(\rho_i^kh_i(p_t^k))-f_{i,t}^k$ and
		$K_i$ as in the first lemma, use $A^{-1}=(a^{ij})$. Multiplication of
		$-\sum_ja^{ij}\Delta w_j=K_ie^{w_i}$ by $\partial_l w_i$ and
		integration by parts give the exact identity
		\begin{align}
			\sum_i\int_{B_r}\partial_l K_i\,e^{w_i}
			=\int_{\partial B_r}\left[\sum_iK_ie^{w_i}\nu_l
			+\sum_{i,j}a^{ij}\left(\partial_\nu w_j\partial_l w_i
			-\tfrac12\nabla w_i\cdot\nabla w_j\nu_l\right)\right].
			\label{I:eq:path-free-translation-Pohozaev}
		\end{align}
		The left side is
		$2\pi\sum_i a_{i,t}^k\partial_l\log K_i(0)+O(E_0^k)$,
		by the disk-mass and first-moment bounds and the uniform $C^2$ norms of
		$\log K_i$.
		
		On the boundary, \eqref{I:eq:path-free-green-C1} writes $w_i$ as
		$-m_{i,t}^{V,k}\log r+R_i(x)$ plus a constant and an $O_{C^1}(E_0^k)$
		error, where $R_i=S_{i,t}^k-f_{i,t}^k$. Its Laplacian is $O(E_0^k)$.
		Replace $R_i$ by its harmonic extension on a slightly larger fixed disk;
		the change in the boundary first derivatives and in $\nabla R_i(0)$ is
		$O(E_0^k)$ by the zero-boundary Poisson estimate. The pure logarithmic
		stress has zero integral, as it is a constant times $\nu_l$.
		The pure harmonic stress has zero divergence and hence zero integral.
		The cross term is exactly
		\[
		-\frac1r\sum_i(A^{-1}m_t^{V,k})_i\partial_l R_i
		=-\frac1r\sum_i a_{i,t}^k\partial_l R_i.
		\]
		The mean-value property for the derivatives of harmonic functions makes
		its integral $-2\pi\sum_i a_{i,t}^k\partial_l R_i(0)$.
		The source boundary term is $O(E_0^k)$ by
		\eqref{I:eq:path-free-exterior}. Thus the exact identity gives
		$\sum_i a_{i,t}^k\partial_l(\log K_i+R_i)(0)=O(E_0^k)$.
		The derivatives of $f_i$ cancel, and the normalized chart has
		$\nabla\widehat\varphi_{p_t^k}(0)=0$. This is
		\eqref{I:eq:path-free-weighted-location}.
		
		Replace every $a_t^k$ by $\sigma$ and every $m_t^{V,k}$ by $m=A\sigma$.
		All coefficients and the required Green derivatives are bounded on the
		fixed separated set, so the error is
		$O(D_k+E_0^k+|\rho^k-\rho^*|)$ by
		\eqref{I:eq:path-free-mass-ledger}. Multiplication by $2\pi N$ and
		\eqref{I:eq:gradient-F} give \eqref{I:eq:path-free-location-estimate}.
	\end{proof}
	
	\begin{theorem}
		\label{I:thm:path-free-coarse-matching}
		Assume (H1)--(H2), (A1), and $h_i\in C^3(M)$ with $h_i>0$.
		Let $u^k$ be
		a sequence of blow-up solutions with $\rho^k\to\rho^*\in\Gamma_N$.
		Suppose either $\rho^k=\rho^*$ or the approach is admissible in
		Definition~\ref{I:def:admissible-path} and
		$\Lambda_{I,N}(\rho^k)>0$. For $N=1$ either transverse side is allowed.
		Choose $\beta$ as above. After a subsequence, the preceding parameters
		satisfy, with $L_k=N^{-1}\sum_t\ell_t^k$,
		$\delta_t^k=\ell_t^k-L_k$, and $E_k=e^{-\beta L_k}$,
		\begin{align}
			|\rho^k-\rho^*|+|\bar a^k-\sigma|&\le CE_k,
			\qquad \max_{s,t}|\ell_t^k-\ell_s^k|\le C,
			\label{I:eq:path-free-mean-comparability}\\
			\max_t|L_kAd_t^k+q\delta_t^k+H_t(p^k)-\overline H(p^k)|
			&\le C(D_k+L_kE_k),
			\label{I:eq:path-free-centered-circle}\\
			|\nabla\ff(p^k)|&\le C(D_k+E_k).
			\label{I:eq:path-free-three-inputs}
		\end{align}
		In particular all three hypotheses of
		Lemma~\ref{I:lem:stationary-profile-matching} are verified, with
		$L_k^2E_k\to0$. Its conclusions give a balanced critical pair
		$(p^*,\delta^*)$ and
		\begin{equation}\label{I:eq:path-free-matching-conclusion}
			\max_t|a_t^k-\sigma|+\max_t|\mu_t^k-\sigma|+|p^k-p^*|\le CE_k,
			\qquad |\delta^k-\delta^*|\le CL_kE_k.
		\end{equation}
		The normalized entire profiles differ from the fixed profile by
		$O_R(E_k)$ in $C^2(B_R)$ for each fixed $R$.
	\end{theorem}
	
	\begin{proof}
		Set $a^k=\rho^k/(2\pi N)$ and $e^k=a^k-\bar a^k$.
		Lemma~\ref{I:lem:path-free-circle} gives $|e^k|\le CE_0^k$.
		For a stationary sequence or the stated lower side, $P(a^k)\ge0$.
		Proposition~\ref{I:prop:radial-variance} therefore yields
		\[
		D_k^2\le CE_0^k,
		\qquad 0\le P(a^k)\le CE_0^k.
		\]
		The second bound follows from the exact identity because its variance
		term is nonpositive. Along an admissible transverse path,
		\[
		D\Lambda_{I,N}(\rho^*)[\dot\rho]
		=-\frac{1}{\pi N}q^T\dot\rho\ne0.
		\]
		Hence $|\rho^k-\rho^*|\le C|P(a^k)|\le CE_0^k$.
		For $N=1$, $D_k=0$ identically and the same exact identity gives
		$|P(a^k)|\le CE_0^k$ on either side. At the stationary parameter the
		parameter estimate is immediate. Thus in every case under consideration
		$|\bar a^k-\sigma|\le CE_0^k$.
		
		The last inequality of
		\eqref{I:eq:path-free-mass-ledger} implies
		\[
		\ell_+^k(D_k+E_0^k)
		\le C(1+\ell_-^k)e^{-\beta\ell_-^k}\longrightarrow0.
		\]
		Subtract
		\eqref{I:eq:path-free-circle-identity} at $s,t$ and use
		$m_{i,t}^{V,k}=m_i+O(D_k+E_0^k)$. Its left scale term becomes
		$q_i(\ell_t^k-\ell_s^k)+o(1)$, whereas all remaining terms are bounded.
		Since $q_i>0$, the logarithmic scale differences are bounded. Consequently
		$E_0^k\le Ce^{-2\beta L_k}\le CE_k$ and $\delta^k$ is bounded.
		
		To extract the centered equation, write
		$m_t^{V,k}=m+A(\bar a^k-\sigma)+Ad_t^k$ in the exact circle identity.
		The regular potential satisfies
		$\log h_i(p_t^k)+S_{i,t}^k(p_t^k)=H_{i,t}(p^k)
		+O(D_k+E_k)$.
		The centered difference of $c_i(a_t^k)$ is $O(D_k)$ by its local
		Lipschitz bound. Centering the scale product gives
		\[
		(m_t^{V,k}-2\mathbf1)\ell_t^k-
		\frac1N\sum_s(m_s^{V,k}-2\mathbf1)\ell_s^k
		=L_kAd_t^k+q\delta_t^k+O(D_k+E_k),
		\]
		because $\delta^k$ is bounded and the mean-mass error is $O(E_k)$.
		These identities give \eqref{I:eq:path-free-centered-circle}.
		Lemma~\ref{I:lem:path-free-location} gives
		\eqref{I:eq:path-free-three-inputs}. Applying
		Lemma~\ref{I:lem:stationary-profile-matching} proves
		\eqref{I:eq:path-free-matching-conclusion}, and the mass identity transfers
		the conclusion to $\mu_t^k$. Smooth dependence of the fixed-component
		entire mass chart, followed by the Lipschitz normalization change, gives
		the final profile assertion.
	\end{proof}
	
	\begin{remark}
		
		If $\widehat m=4$, the preceding theorem gives
		$|m_{i,t}^{V,k}-4|\le Ce^{-\beta L_k}\le C/L_k$ for every component and
		bubble. Thus the mass-rate hypothesis of the local
		\cite[Theorem~4.3, (4.55)--(4.56)]{LinZhang2013} is now available at
		stationary parameters as well. In either regime, for
		$1\le R\le r_0/\epsilon_t^k$,
		\[
		R^{-(m_{i,t}^{V,k}-m_i)}=1+O(L_ke^{-\beta L_k}).
		\]
		Indeed $|(m_{i,t}^{V,k}-m_i)\log R|\le CL_ke^{-\beta L_k}\to0$;
		the mean-value theorem for the exponential proves the displayed formula.
		These facts permit the second-order local expansion in
		Section~\ref{I:sec:integrated-critical-estimates}, including the uniform
		integrations up to a fixed physical radius.
	\end{remark}
	
	For the preliminary parameters, define
	\begin{equation}\label{I:eq:def-preliminary-residual}
		\begin{aligned}
			\widehat R_i^k:={}&\rho_i^kh_ie^{\Theta_i^k}
			-\sum_t\chi_{\widehat p_t^k}e^{-\varphi_{\widehat p_t^k}}
			e^{W_{i,\widehat p_t^k,\widehat\lambda_t^k}}-\rho_i^k\\
			&+\sum_t\int_M\chi_{\widehat p_t^k}e^{-\varphi_{\widehat p_t^k}}
			e^{W_{i,\widehat p_t^k,\widehat\lambda_t^k}}\,dV_g.
		\end{aligned}
	\end{equation}
	Then $\int_M\widehat R_i^k\,dV_g=0$.
	
	\begin{proposition}
		\label{I:prop:coarse-source-entry}
		Under the hypotheses of Theorem~\ref{I:thm:path-free-coarse-matching},
		take $\widehat p_t^k=p_t^k$, $\widehat\lambda_t^k=e^{\ell_t^k}$ in
		\eqref{I:eq:def-preliminary-residual}. If $1<p<1+\beta/2$ and
		$\kappa_p=(\beta-2(p-1))/p>0$, then
		\begin{equation*}
			\begin{aligned}
				\sum_i\|\widehat R_i^k\|_{L^1}\le Ce^{-\beta L_k},\qquad
				&\sum_i\|\widehat R_i^k\|_{L^p}\le Ce^{-\kappa_pL_k},\\
				\|u^k-U_{\widehat p^k,\widehat\lambda^k,0}\|_{W^{2,p}(M)^n}
				&\le Ce^{-\kappa_pL_k}.
			\end{aligned}
		\end{equation*}
	\end{proposition}
	
	\begin{proof}
		The local mass calculation proves more than equality of integrals: it
		bounds the $L^1$ difference of the two positive source densities on
		each disk by $C(\epsilon_t^k)^{2\beta}$. The normalized entire family
		also satisfies, uniformly near $\sigma$,
		\[
		\int_{\R^2}|e^{V_i(a)}-e^{v_i}|\,dy\le C|a-\sigma|.
		\]
		This follows from Lemma~\ref{I:lem:radial-family-derivatives} and
		the locally Lipschitz maximum-zero normalization.
		Theorem~\ref{I:thm:path-free-coarse-matching} makes this error $O(E_k)$.
		The fixed cutoff removes a tail of size $O(E_k)$, and
		\eqref{I:eq:path-free-exterior} controls the actual exterior source.
		The factor $e^{-\varphi_{p_t^k}}$ in the projected density cancels the
		volume Jacobian exactly. Thus the uncentered residual has $L^1$ norm
		$O(E_k)$. Subtracting its mean at most doubles this norm.
		
		The same local representation and boundedness of the normalized entire
		profiles give $\|\widehat R_i^k\|_\infty\le Ce^{2L_k}$, using the
		proved scale comparison. Therefore interpolation gives
		\[
		\|\widehat R_i^k\|_p
		\le \|\widehat R_i^k\|_1^{1/p}
		\|\widehat R_i^k\|_\infty^{1-1/p}
		\le Ce^{-[\beta-2(p-1)]L_k/p}.
		\]
		Finally both $u_i^k$ and the projected bubbles have mean zero, and their
		defining equations give exactly
		$-\Delta_g(u_i^k-U_{i,\widehat p^k,\widehat\lambda^k,0})
		=\sum_j a_{ij}\widehat R_j^k$.
		The zero-mean elliptic estimate proves the last bound. Since $p>1$ on
		a two-dimensional surface, it also proves strong convergence in $\mathcal H$.
	\end{proof}
	
	\subsection{Coefficients in the mass expansion}
	
	Choose a measurable partition $\{\Omega_t\}_{t=1}^N$ of $M$, up to
	boundaries of zero area, with each $\Omega_t$ containing a fixed disk
	about $p_t$ and no other concentration point in its closure. For nearby
	configurations, keep the boundaries fixed away from these disks. Set
	\[
	\mathcal G(x;p):=\sum_sG^g(x,p_s),\qquad
	\mathcal G_t^*(p):=R^g(p_t)+\sum_{s\ne t}G^g(p_t,p_s).
	\]
	For $2<m_i<4$, define
	\begin{equation}\label{I:def:Dit}
		D_{i,t}(p):=\lim_{r\downarrow0}\left[
		-\frac{2\pi}{m_i-2}r^{2-m_i}
		+\int_{\Omega_t\setminus U_r(p_t)}
		\frac{h_i(x)}{h_i(p_t)}
		e^{2\pi m_i(\mathcal G(x;p)-\mathcal G_t^*(p))}\,dV_g
		\right].
	\end{equation}
	In the normalized chart, the integrand is $|y|^{-m_i}B(y)\,dy$, where
	$B(0)=1$ and $B(y)=1+DB(0)y+O(|y|^2)$. The linear term has zero
	integral on centered annuli, and
	\[
	\int_{|y|<r_0}|y|^{-m_i}|B(y)-1-DB(0)y|\,dy
	\le C\int_0^{r_0}r^{3-m_i}\,dr<\infty.
	\]
	Thus $D_{i,t}$ exists and is continuous on separated configurations
	with locally fixed partition boundaries.
	
	When $\widehat m=4$, define
	\begin{equation}\label{I:def:Lit}
		L_{i,t}^{\mathrm{corr}}(p):=
		\frac{e^{\fJ_i}}4\left[
		\frac{\Delta_{g,x}\rf_{p,t}^i(p_t)}{\rf_{p,t}^i(p_t)}-2K_g(p_t)
		\right].
	\end{equation}
	For $N=1$, this is the coefficient in
	\cite[Theorem~1.3]{LinZhang2013}.
	
	\begin{lemma}\label{I:lem:radial-second-moment}
		If $m_i=4$, then
		\[
		\int_{|y|<R}|y|^2e^{v_i(y)}\,dy
		=2\pi e^{\fJ_i}\log R+O(1).
		\]
		For $A_{i,t}(z):=e^{\widehat\varphi_{p_t}(z)}
		\rf_{p,t}^i(y_{p_t}^{-1}(z))/\rf_{p,t}^i(p_t)$,
		\[
		\frac{e^{\fJ_i}}4\Delta A_{i,t}(0)=L_{i,t}^{\mathrm{corr}}(p).
		\]
	\end{lemma}
	
	\begin{proof}
		The tail $e^{v_i(r)}=e^{\fJ_i}r^{-4}(1+O(r^{-\eta}))$ gives the
		moment formula by polar integration. Since
		$\widehat\varphi_{p_t}(0)=0$, $\nabla\widehat\varphi_{p_t}(0)=0$,
		and $\Delta\widehat\varphi_{p_t}(0)=-2K_g(p_t)$,
		\[
		\Delta A_{i,t}(0)
		=\frac{\Delta_{g,x}\rf_{p,t}^i(p_t)}{\rf_{p,t}^i(p_t)}-2K_g(p_t).
		\]
		The angular average of the quadratic Taylor term is
		$\frac14\Delta A_{i,t}(0)|z|^2$, which gives the coefficient
		in \eqref{I:def:Lit}.
	\end{proof}
	
	Equivalently, \eqref{I:def:tilde-G} gives
	$\Delta_g\widetilde G_t=N$ near $p_t$, and hence
	\[
	\begin{aligned}
		L_{i,t}^{\mathrm{corr}}(p)=e^{\fJ_i}\Bigl[
		&\frac{\Delta_gh_i(p_t)}{4h_i(p_t)}-\frac12K_g(p_t)+2\pi N\\
		&+4\pi\nabla_g\log h_i(p_t)\cdot\nabla_x\widetilde G_t(p_t)
		+16\pi^2|\nabla_x\widetilde G_t(p_t)|_g^2\Bigr].
	\end{aligned}
	\]
	
	Set
	\begin{equation}\label{I:def:tail-coefficients}
		\mathcal L_t(p):=
		\begin{cases}
			\displaystyle\frac{\widehat m-2}{\pi N}
			\sum_{i\in\widehat I}\sum_s e^{\fJ_i}
			\frac{\rf_{p,s}^i(p_s)}{\rf_{p,t}^i(p_t)}D_{i,s}(p),
			&\widehat m<4,\\[2mm]
			\displaystyle\frac4N\sum_{i,s}
			\frac{\rf_{p,s}^i(p_s)}{\rf_{p,t}^i(p_t)}L_{i,s}^{\mathrm{corr}}(p),
			&\widehat m=4.
		\end{cases}
	\end{equation}
	For a balanced pair $(p^*,r^*)$, define
	\begin{equation}\label{I:def:L-balanced}
		\mathcal L_\infty(p^*,r^*):=
		e^{-(\widehat m-2)\delta_t^*}\mathcal L_t(p^*).
	\end{equation}
	
	\begin{proposition}\label{I:prop:balanced-coefficient-compatibility}
		For every balanced pair, $\mathcal L_\infty$ is independent of $t$
		and of the partition. Moreover,
		\begin{equation}\label{I:eq:L-compatible}
			e^{-(\widehat m-2)\delta_t^*}\mathcal L_t(p^*)
			=e^{-(\widehat m-2)\delta_s^*}\mathcal L_s(p^*),
			\qquad s,t\in I'.
		\end{equation}
	\end{proposition}
	
	\begin{proof}
		By \eqref{I:eq:balance-weights},
		\[
		\frac{\rf_{p^*,s}^i(p_s^*)}{\rf_{p^*,t}^i(p_t^*)}
		=e^{-q_i(\delta_s^*-\delta_t^*)}.
		\]
		Substitution gives
		\[
		\mathcal L_\infty=
		\begin{cases}
			\displaystyle\frac{\widehat m-2}{\pi N}
			\sum_{i\in\widehat I,s}e^{\fJ_i}
			e^{-(\widehat m-2)\delta_s^*}D_{i,s}(p^*),&\widehat m<4,\\[2mm]
			\displaystyle\frac4N\sum_{i,s}e^{-2\delta_s^*}
			L_{i,s}^{\mathrm{corr}}(p^*),&\widehat m=4,
		\end{cases}
		\]
		which proves \eqref{I:eq:L-compatible}. Below four, the density
		in each weighted finite part is
		\[
		\frac{e^{-q_i\delta_t^*}}{\rf_{p^*,t}^i(p_t^*)}
		h_i(x)e^{2\pi m_i\mathcal G(x;p^*)},\qquad i\in\widehat I.
		\]
		Its prefactor is independent of $t$. Reassigning a region away from the
		concentration points therefore preserves the weighted sum; the local
		subtractions are unchanged. At four the coefficient contains no partition.
	\end{proof}
	
	\subsection{Integrated mass and moment estimates}
	\label{I:sec:integrated-critical-estimates}
	
	In this subsection $p_t^k=\widehat p_t^k$ and
	$\epsilon_t^k=\varepsilon_{k,t}$ are the preliminary local parameters;
	$L_k=N^{-1}\sum_t\log(\epsilon_t^k)^{-1}$. The exact minimizing
	parameters are introduced in Section~\ref{I:sec:entry-theorem}.
	
	The entire masses, partition masses and disk masses are, respectively,
	\begin{equation*}
		a_{i,t}^k=\frac1{2\pi}\int_{\R^2}e^{V_i^{k,t}},\qquad
		\mu_{i,t}^k=\frac1{2\pi}\int_{\Omega_t^k}s_i^k\,dV_g,\qquad
		\nu_{i,t}^k(r)=\frac1{2\pi}\int_{U_r(p_t^k)}s_i^k\,dV_g.
	\end{equation*}
	Thus
	\[
	P(a_t^k)=0,\qquad \sum_t\mu_t^k=\frac{\rho^k}{2\pi},
	\qquad m_{i,t}^{V,k}=(Aa_t^k)_i.
	\]
	Proposition~\ref{I:prop:refined-local-Harnack} uses the same radial
	profiles as the first-order approximation, with a displacement and an
	angular correction. Their masses are unchanged. Scale comparability gives
	$\Psi_{\widehat m}((\epsilon_t^k)^{-1})\asymp\Psi_k$.
	
	\begin{lemma}
		
		In the chart at $p_t^k$ set
		\[
		J_{i,t}^k(x)=\frac{h_i(x)}{h_i(p_t^k)}
		e^{\widehat\varphi_{p_t^k}(x)+f_{i,t}^k(x)+\phi_{i,t}^k(x)}.
		\]
		Here and below a function on the surface is composed with the inverse
		chart when its argument is a coordinate vector. The functions $J_{i,t}^k$
		and their reciprocals have uniform $C^3$ bounds on a smaller fixed disk,
		and $J_{i,t}^k(0)=1$. On still smaller disks, after identifying the moving
		charts with the limiting chart, they converge in $C^2$ to
		\begin{equation}\label{I:eq:limiting-J}
			B_{i,t}(x)=e^{\widehat\varphi_{p_t^*}(x)}
			\frac{\rf_{p^*,t}^i(x)}{\rf_{p^*,t}^i(p_t^*)}.
		\end{equation}
		In the translated scaled coordinate $z=x/\epsilon_t^k-\zeta_t^k$,
		$|\zeta_t^k|\le C\epsilon_t^k$, the actual source measure has density
		\begin{equation}\label{I:eq:common-density}
			J_{i,t}^k(\epsilon_t^k(z+\zeta_t^k))
			e^{V_i^{k,t}(z)+\Phi_{i,t}^k(z)+R_{i,t}^k(z)}\,dz.
		\end{equation}
		The angular mean of $\Phi_{i,t}^k$ is zero on every centered circle.
		For $q_*<2$ the bounds are
		\begin{equation*}
			|\Phi_{i,t}^k(z)|\le C\epsilon_t^k r(1+r)^{-q_*},\qquad
			|R_{i,t}^k(z)|\le C(\epsilon_t^k)^2
			(1+r)^{2-q_*}\log(2+r),\quad r=|z|.
		\end{equation*}
		For $q_*=2$ they are
		\begin{equation}\label{I:eq:integrated-pointwise-critical}
			|\Phi_{i,t}^k(z)|\le C\epsilon_t^k r(1+r)^{-2},\qquad
			|R_{i,t}^k(z)|\le C(\epsilon_t^k)^2\log^2(2+r).
		\end{equation}
		The radial densities satisfy
		\begin{equation}\label{I:eq:uniform-sharp-density-tail}
			e^{V_i^{k,t}(z)}\le C(1+|z|)^{-m_i},\qquad
			e^{V_i^{k,t}(z)}=e^{c_i^k}|z|^{-m_{i,t}^{V,k}}
			(1+O(|z|^{-q_*/2}))\quad (|z|\ge1),
		\end{equation}
		where $c_i^k\to\fJ_i$. The first bound is uniform for
		$|z|\le c/\epsilon_t^k$, and the tail expansion holds uniformly for all
		$|z|\ge1$.
	\end{lemma}
	\begin{proof}
		The exact change of variables in the source equation gives
		\eqref{I:eq:common-density}. Proposition~\ref{I:prop:refined-local-Harnack}
		gives the displacement, angular mean and remainder bounds, initially with
		$q_*^k=\min_i m_{i,t}^{V,k}-2$ in place of $q_*$. The coarse theorem gives
		$|m_{i,t}^{V,k}-m_i|\le Ce^{-\beta L_k}$; consequently
		\[
		\sup_{1\le r\le c/\epsilon_t^k}
		|(m_{i,t}^{V,k}-m_i)\log r|\le CL_ke^{-\beta L_k}\longrightarrow0.
		\]
		This proves the stated bounds without a loss in the limiting exponent.
		The uniform radial tail expansion gives
		\eqref{I:eq:uniform-sharp-density-tail}. The displacement is retained
		in the argument of $J_{i,t}^k$; the radial normalization is unchanged.
		
		It remains to identify $J$. The coarse Green formula on the boundary
		of a fixed disk gives the boundary oscillation of
		$\Theta_i^k-f_{i,t}^k$ as that of $S_{i,t}^k-f_{i,t}^k$ plus
		$O(e^{-\beta L_k})$. Moreover
		\[
		\Delta(S_{i,t}^k-f_{i,t}^k)
		=\left(2\pi\sum_s m_{i,s}^{V,k}-(A\rho^k)_i\right)
		e^{\widehat\varphi_{p_t^k}}=O(e^{-\beta L_k}).
		\]
		Subtract its zero-boundary Poisson solution. The remaining function is
		harmonic, and its zero-center normalization agrees with that of
		$\phi_{i,t}^k$. The maximum principle and interior derivative estimates
		for harmonic functions give, on a smaller disk,
		\[
		f_{i,t}^k+\phi_{i,t}^k
		=S_{i,t}^k-S_{i,t}^k(0)+O_{C^2}(e^{-\beta L_k}).
		\]
		The Poisson solution has this $C^2$ bound by its smooth right side and
		fixed disk. Since $S_{i,t}^k\to2\pi m_i\widetilde G_t(\cdot;p^*)$
		smoothly there, \eqref{I:eq:limiting-J} follows. Uniform $C^3$ bounds
		on smaller disks follow from the smooth background, the uniform
		boundary oscillation and interior harmonic derivative estimates.
	\end{proof}
	
	\begin{lemma}
		\label{I:lem:integrated-angular-moments}
		For every sufficiently small fixed $r>0$,
		\begin{equation}\label{I:eq:critical-moment-bounds}
			\left|\int_{U_r(p_t^k)}y_{p_t^k}(x)s_i^k\,dV_g\right|
			+\int_{U_r(p_t^k)}|y_{p_t^k}(x)|^2s_i^k\,dV_g
			\le C_r\Psi_k.
		\end{equation}
		In the unweighted local mass integral, replacing the density in
		\eqref{I:eq:common-density} by
		$J_{i,t}^k(\epsilon_t^kz)e^{V_i^{k,t}(z)}$ introduces
		$o((\epsilon_t^k)^{q_i})$ if $q_*<2$ and $q_i<2$,
		$o(\Psi_k)$ if $q_*<2$ and $q_i\ge2$, and
		$O((\epsilon_t^k)^2)=o(\Psi_k)$ if $q_*=2$.
		The integration disk may be centered at either of the two centers.
	\end{lemma}
	\begin{proof}
		Suppress $k,t$, set  $T=r/\epsilon$ and $L_\epsilon=\log\epsilon^{-1}$.
		The pointwise bounds imply $\|\Phi\|_\infty+\|R\|_\infty=o(1)$ on
		this disk. Thus
		$|e^{\Phi+R}-1-\Phi|\le C(|R|+\Phi^2)$.
		The angular integral of $e^{V_i}\Phi_i$ is exactly zero. The remaining
		unweighted errors are bounded by
		\begin{equation*}
			C\int_{|z|<T}e^{V_i}
			\bigl(|R_i|+\Phi_i^2+\epsilon|z|\,|\Phi_i|\bigr)\,dz.
		\end{equation*}
		For $q_*<2$, each of the first and third terms is at most
		\[
		C\epsilon^2\int_0^T
		(1+s)^{2-q_*-m_i}\log(2+s)\,s\,ds
		\le C_r\left[\epsilon^2(1+L_\epsilon^2)
		+\epsilon^{q_*+q_i}(1+L_\epsilon^2)\right].
		\]
		For the squared angular term replace $q_*$ by $2q_*$ in the last bound.
		To verify these inequalities, split the integral at $1$; on $[1,T]$
		its power is $s^{1-q_*-q_i}$, and integrate it, with a logarithmic
		factor when $q_*+q_i=2$. The strict inequalities $q_*>0$, $q_i>0$
		show that these bounds are $o(\epsilon^{q_i})$ when $q_i<2$, and
		$o(\epsilon^{q_*})$ when $q_i\ge2$. At $q_*=2$ all three integrals
		are $O(\epsilon^2)$, since
		$\int_0^\infty(1+s)^{-4}\log^2(2+s)s\,ds<\infty$.
		
		The displacement changes $J$ by $O(\epsilon^2)$ and the physical
		center by $O(\epsilon^2)$. The symmetric difference of the two disks
		has area $O_r(\epsilon^2)$ in physical coordinates. On its boundary
		the actual and radial physical densities are $O_r(\epsilon^{q_i})$,
		by \eqref{I:eq:uniform-sharp-density-tail} and the uniform smallness
		of $\Phi,R$. Its unweighted mass is therefore
		$O_r(\epsilon^{q_i+2})$. These errors have precisely the orders claimed.
		
		For the second moment, boundedness of $J$ and $e^{\Phi+R}$ gives
		\[
		\int_{U_r}|x|^2s_i\,dV_g
		\le C_r\epsilon^2\int_0^T s^3(1+s)^{-m_i}ds+C\epsilon^4.
		\]
		This is $O_r(\epsilon^{q_i})$ if $q_i<2$,
		$O_r(\epsilon^2(1+L_\epsilon))$ if $q_i=2$, and
		$O_r(\epsilon^2)$ if $q_i>2$. Each is $O_r(\Psi_k)$.
		For the signed first moment the radial constant term vanishes. The
		$J-1$ term is bounded by this same second-moment integral. The terms
		from $\Phi$ and $R+\Phi^2$ are bounded, respectively, by
		\[
		C\epsilon^2\int_0^T s^3(1+s)^{-m_i-q_*}ds,
		\quad C\epsilon^3\int_0^T
		s^2(1+s)^{2-q_*-m_i}\log(2+s)ds
		\]
		and, for the squared angular term, by
		$C\epsilon^3\int_0^T s^4(1+s)^{-m_i-2q_*}ds$.
		Integrating the powers
		shows that these are bounded by
		$C_r[\epsilon^2(1+L_\epsilon^3)
		+\epsilon^{q_*+q_i}(1+L_\epsilon^3)]$ when $q_*<2$; this is
		$o(\epsilon^{q_*})$. At $q_*=2$ the bounds in
		\eqref{I:eq:integrated-pointwise-critical} instead give $O(\epsilon^2)$.
		The physical displacement adds $O(\epsilon^2)$ times the bounded
		total mass. This proves \eqref{I:eq:critical-moment-bounds}.
	\end{proof}
	
	\begin{proposition}
		\label{I:prop:integrated-mass-ledger}
		Under the sequence hypotheses of Theorem~\ref{I:thm:path-free-coarse-matching},
		set  $\Delta_{i,t}^k=\mu_{i,t}^k-a_{i,t}^k$ for a fixed local partition,
		transported with the charts near the centers, and denote its limiting
		partition at $p^*$ by $\{\Omega_t^*\}$. If $q_*<2$, then
		\begin{align}
			\Delta_{i,t}^k
			&=\frac{e^{\fJ_i}}{2\pi}D_{i,t}(p^*)
			(\epsilon_t^k)^{q_i}+o((\epsilon_t^k)^{q_i}),&&q_i<2,
			\label{I:eq:sharp-integrated-nonlocal}\\
			\Delta_{i,t}^k&=o(\Psi_k),&&q_i\ge2.
			\label{I:eq:sharp-integrated-higher}
		\end{align}
		If $q_*=2$, then
		\begin{equation}\label{I:eq:sharp-integrated-borderline}
			\Delta_{i,t}^k=L_{i,t}^{\mathrm{corr}}(p^*)
			(\epsilon_t^k)^2\log(\epsilon_t^k)^{-1}+o(\Psi_k).
		\end{equation}
		In particular $|\Delta_{i,t}^k|\le C\Psi_k$ in both regimes, and
		$|\nu_{i,t}^k(r)-a_{i,t}^k|\le C_r\Psi_k$ for every fixed small $r$.
	\end{proposition}
	\begin{proof}
		For fixed $r$, the exact decomposition is
		\begin{equation}\label{I:eq:exact-core-tail-exterior}
			\begin{split}
				2\pi\Delta_{i,t}^k={}&
				\int_{|z|<r/\epsilon_t^k}
				(J_{i,t}^k(\epsilon_t^kz)-1)e^{V_i^{k,t}(z)}dz
				-\int_{|z|>r/\epsilon_t^k}e^{V_i^{k,t}(z)}dz\\
				&+\int_{\Omega_t^k\setminus U_r(p_t^k)}s_i^k\,dV_g
				+\mathcal E_{i,t}^k(r),
			\end{split}
		\end{equation}
		where Lemma~\ref{I:lem:integrated-angular-moments} gives the stated
		little-oh orders for $\mathcal E$, and $O((\epsilon_t^k)^2)$ at
		$q_*=2$. The entire mass has not been replaced by its limiting value.
		
		Suppose $q_i<2$. Subtract the linear Taylor term of $J^k-1$ inside
		the radial integral; its integral is zero. In physical coordinates,
		\eqref{I:eq:uniform-sharp-density-tail} and coarse exponent control give
		\[
		(\epsilon_t^k)^{-q_i-2}e^{V_i^{k,t}(x/\epsilon_t^k)}
		\longrightarrow e^{\fJ_i}|x|^{-m_i}\quad(x\ne0).
		\]
		After multiplication by the Taylor remainder, the integrand is bounded
		by $C|x|^{2-m_i}$ outside $|x|<\epsilon_t^k$. The normalized integral
		on that inner disk is $O((\epsilon_t^k)^{2-q_i})$. Hence dominated
		convergence, using $m_i<4$, gives
		\[
		(\epsilon_t^k)^{-q_i}\int_{|z|<r/\epsilon_t^k}
		(J_{i,t}^k(\epsilon_t^kz)-1)e^{V_i^{k,t}(z)}dz
		\longrightarrow e^{\fJ_i}\int_{|x|<r}
		\frac{B_{i,t}(x)-1-DB_{i,t}(0)x}{|x|^{m_i}}dx.
		\]
		The tail term divided by $(\epsilon_t^k)^{q_i}$ tends to
		$2\pi e^{\fJ_i}r^{-q_i}/q_i$.
		
		For the exterior use the coarse circle identity at bubble $t$ and
		the coarse Green representation away from every center. They give,
		locally uniformly on that exterior,
		\begin{equation*}
			\frac{s_i^k(x)}{(\epsilon_t^k)^{q_i}}
			\longrightarrow e^{\fJ_i}\frac{h_i(x)}{h_i(p_t^*)}
			e^{2\pi m_i(\mathcal G(x;p^*)-\mathcal G_t^*(p^*))}.
		\end{equation*}
		Indeed the logarithm of the quotient has additional term
		$-(m_{i,t}^{V,k}-m_i)\ell_t^k=o(1)$, which is exactly why the coarse
		mass rate was established before this calculation. On each fixed
		exterior the convergence is uniform and supplies an integrable bound.
		By \eqref{I:def:Dit} and the vanishing angular integral of the linear
		Taylor term, for every fixed small $r$,
		\[
		\begin{split}
			D_{i,t}(p^*)={}&-\frac{2\pi}{q_i}r^{-q_i}
			+\int_{\Omega_t^*\setminus U_r(p_t^*)}
			\frac{h_i(x)}{h_i(p_t^*)}
			e^{2\pi m_i(\mathcal G(x;p^*)-\mathcal G_t^*(p^*))}\,dV_g\\
			&+\int_{|x|<r}
			\frac{B_{i,t}(x)-1-DB_{i,t}(0)x}{|x|^{m_i}}dx.
		\end{split}
		\]
		The core, terminal tail, and exterior limits therefore sum to
		$e^{\fJ_i}D_{i,t}(p^*)$, proving
		\eqref{I:eq:sharp-integrated-nonlocal}.
		
		If $q_*<2$ and $q_i\ge2$, the core Taylor integral is
		$O_r((\epsilon_t^k)^2(1+\log(\epsilon_t^k)^{-1}))$ by the second
		moment. The entire tail and exterior are $O_r((\epsilon_t^k)^{q_i})$
		by the same radial and circle formulas. All are $o(\Psi_k)$, as are
		the errors in the preceding lemma. This proves
		\eqref{I:eq:sharp-integrated-higher}.
		
		Finally suppose $q_*=2$, so every $m_i=4$. The entire tail and the
		fixed exterior are $O_r((\epsilon_t^k)^2)=o(\Psi_k)$. Taylor-expand
		$J^k$ to second order in the core. Its linear integral is zero, and
		the third-order remainder has integral $O_r((\epsilon_t^k)^2)$:
		in polar coordinates it is bounded by
		$C(\epsilon_t^k)^3\int_0^{r/\epsilon_t^k}s^4(1+s)^{-4}ds$.
		Moreover coarse exponent control and the radial expansion give
		\[
		\int_{|z|<r/\epsilon_t^k}|z|^2e^{V_i^{k,t}(z)}dz
		=2\pi e^{\fJ_i}\log(\epsilon_t^k)^{-1}(1+o(1)).
		\]
		To check this last assertion, integrate the leading tail on
		$1<s<r/\epsilon_t^k$ and use
		$\sup|s^{-(m_{i,t}^{V,k}-4)}-1|=o(1)$; the tail remainder has a
		bounded second moment. The quadratic integral in
		\eqref{I:eq:exact-core-tail-exterior}, divided by $2\pi$, is therefore
		\[
		\frac14 e^{\fJ_i}\Delta B_{i,t}(0)
		(\epsilon_t^k)^2\log(\epsilon_t^k)^{-1}+o(\Psi_k).
		\]
		Lemma~\ref{I:lem:radial-second-moment} identifies this coefficient
		with $L_{i,t}^{\mathrm{corr}}(p^*)$. This proves
		\eqref{I:eq:sharp-integrated-borderline}. Omitting the fixed exterior
		in the exact decomposition proves the disk bound in all cases.
	\end{proof}

	\subsection{Sharp matching and the parameter identity}
	
	\begin{proposition}
		\label{I:prop:sharp-bootstrap-closure}
		Under the sequence hypotheses of Theorem~\ref{I:thm:path-free-coarse-matching},
		the three inputs of Lemma~\ref{I:lem:stationary-profile-matching} hold
		with $E_k=\Psi_k$. In particular,
		\begin{align}
			|\rho^k-\rho^*|+\max_t|a_t^k-\sigma|
			+\max_t|\mu_t^k-\sigma|+|p^k-p^*|&\le C\Psi_k,
			\label{I:eq:sharp-closed-masses}\\
			|\delta^k-\delta^*|&\le CL_k\Psi_k.
			\label{I:eq:sharp-closed-scales}
		\end{align}
		On every fixed matching annulus,
		\begin{equation}\label{I:eq:sharp-green-C1}
			\Theta_i^k=\overline\Theta_i^k+
			2\pi\sum_s m_{i,s}^{V,k}G^g(\cdot,p_s^k)+O_{C^1,r}(\Psi_k).
		\end{equation}
		For every fixed $R$, the maximum-zero entire profiles satisfy
		$\|V^{k,t}-v\|_{C^2(B_R)}\le C_R\Psi_k$.
		With $\rf_{p,t}^i$ defined in \eqref{I:def:rf_i}, the relative scales satisfy
		\begin{equation}\label{I:eq:coefficient-balance}
			\frac{\rf_{p^k,t}^i(p_t^k)}{\rf_{p^k,s}^i(p_s^k)}
			=\left(\frac{\epsilon_t^k}{\epsilon_s^k}\right)^{m_i-2}
			\bigl(1+O(L_k\Psi_k)\bigr),
			\qquad i\in I,\ s,t\in I'.
		\end{equation}
	\end{proposition}
	\begin{proof}
		Proposition~\ref{I:prop:integrated-mass-ledger} gives
		$e^k=\rho^k/(2\pi N)-\bar a^k=N^{-1}\sum_t\Delta_t^k=O(\Psi_k)$.
		On the stationary parameter or the lower side, the nonpositive variance
		in Proposition~\ref{I:prop:radial-variance} gives
		$D_k^2\le C\Psi_k$ and $0\le P(\rho^k/(2\pi N))\le C\Psi_k$.
		Transversality then gives $|\rho^k-\rho^*|\le C\Psi_k$. For $N=1$
		the variance is zero and the same exact formula bounds the absolute
		wall defect on either side. Thus $|\bar a^k-\sigma|\le C\Psi_k$.
		
		To prove \eqref{I:eq:sharp-green-C1}, Taylor-expand $G^g(x,\cdot)$
		and its first $x$ derivatives to second order on core disks disjoint
		from the matching annulus. Their linear integrals use the signed first
		moment in \eqref{I:eq:critical-moment-bounds}; the remainders use its
		absolute second moment. Replacing each core mass by $2\pi a_s^k$
		costs $O_r(\Psi_k)$ by the disk bound. On the remaining region,
		the circle formula and coarse exponent control give
		$\|s_i^k\|_\infty\le C_r(\epsilon_1^k)^{q_i}\le C_r\Psi_k$.
		The Green kernel and its first derivatives have uniformly bounded
		absolute integrals there. This proves the sharp Green formula, including
		the differentiated remainder.
		
		Average the refined local formula on a fixed circle. The angular term
		has zero mean at the shifted center; changing the circle center by
		$O((\epsilon_t^k)^2)$ costs $O_r((\epsilon_t^k)^2)$, using the
		first derivative bounds and the radial tail. The remainder has boundary
		size $O_r(L_k\Psi_k)$: below four it is
		$O_r((\epsilon_t^k)^{q_*}\log(\epsilon_t^k)^{-1})$, and at four it
		is $O_r((\epsilon_t^k)^2\log^2(\epsilon_t^k)^{-1})$.
		The mean-value calculation for $S_{i,t}^k-f_{i,t}^k$ now has error
		$O_r(\Psi_k)$. Hence the circle identity
		\eqref{I:eq:path-free-circle-identity} holds with error
		$O_r(L_k\Psi_k)$. Center it in $t$. Lipschitz dependence of $c_i(a_t^k)$,
		the bounded relative scales and $|\bar a^k-\sigma|=O(\Psi_k)$ give
		\[
		\max_t|L_kAd_t^k+q\delta_t^k+H_t(p^k)-\overline H(p^k)|
		\le C(D_k+L_k\Psi_k).
		\]
		
		In the exact translation identity
		\eqref{I:eq:path-free-translation-Pohozaev}, Taylor-expand
		$\partial_l\log K_i$ to first order. Its linear term uses the
		signed first moment and its error uses the second moment, since
		$\log K_i$ has uniform $C^3$ bounds. The left side is therefore
		$2\pi\sum_i a_{i,t}^k\partial_l\log K_i(0)+O_r(\Psi_k)$.
		Use \eqref{I:eq:sharp-green-C1} on the boundary. Its regular part
		$S_{i,t}^k-f_{i,t}^k$ has Laplacian $O(\Psi_k)$ by the mass identity;
		replace it by its harmonic extension as in
		Lemma~\ref{I:lem:path-free-location}. The pure logarithmic and harmonic
		stresses have zero integral. Their cross stress gives exactly
		$-2\pi\sum_i a_{i,t}^k\nabla(S_{i,t}^k-f_{i,t}^k)(0)$.
		The boundary source is $O_r(\Psi_k)$. Cancellation of $f_i$ proves
		\[
		\sum_i a_{i,t}^k[\nabla\log h_i(p_t^k)+\nabla S_{i,t}^k(p_t^k)]
		=O_r(\Psi_k),\qquad
		|\nabla\ff(p^k)|\le C(D_k+\Psi_k).
		\]
		Since $L_k^2\Psi_k\to0$, Lemma~\ref{I:lem:stationary-profile-matching}
		gives \eqref{I:eq:sharp-closed-masses}--\eqref{I:eq:sharp-closed-scales}.
		The mass chart gives the fixed-ball profile estimate, and the centered
		circle equation now gives the balance error $O(L_k\Psi_k)$.
	\end{proof}
	
	\begin{theorem}
		\label{I:thm:sharp-global-mass-identity}
		For every sequence in Theorem~\ref{I:thm:path-free-coarse-matching},
		\begin{equation}\label{I:eq:global-identity-from-ledger}
			\Lambda_{I,N}(\rho^k)=-\frac2N\sum_{i,t}q_i\Delta_{i,t}^k
			+O(\Psi_k^2).
		\end{equation}
		Consequently, for each $t$,
		\begin{equation}\label{I:eq:sharp-global-preliminary}
			\Lambda_{I,N}(\rho^k)
			=-\mathcal L_t(p^*)\Psi_{\widehat m}(\widehat\lambda_t^k)+o(\Psi_k).
		\end{equation}
		The balanced version uses $\mathcal L_\infty$ from
		\eqref{I:def:L-balanced} and the mean logarithmic preliminary scale.
	\end{theorem}
	\begin{proof}
		Set $a^k=\rho^k/(2\pi N)$ and $e^k=a^k-\bar a^k$.
		By \eqref{I:def:Lam_rho}, $P(a^k)=\Lambda_{I,N}(\rho^k)$.
		The exact variance identity reads
		\[
		P(a^k)=\frac1N\sum_t(d_t^k)^TAd_t^k-2q^Te^k
		-2[A(\bar a^k-\sigma)]^Te^k-(e^k)^TAe^k.
		\]
		By Proposition~\ref{I:prop:sharp-bootstrap-closure}, every term except
		$-2q^Te^k$ is $O(\Psi_k^2)$. The exact identity
		$e^k=N^{-1}\sum_t\Delta_t^k$ proves
		\eqref{I:eq:global-identity-from-ledger}.
		For $q_*<2$, insert \eqref{I:eq:sharp-integrated-nonlocal}--
		\eqref{I:eq:sharp-integrated-higher}; only $i\in\widehat I$ survive.
		The result is
		\[
		-\frac{q_*}{\pi N}\sum_{i\in\widehat I,t}
		e^{\fJ_i}D_{i,t}(p^*)(\epsilon_t^k)^{q_*}+o(\Psi_k).
		\]
		For $q_*=2$, \eqref{I:eq:sharp-integrated-borderline} instead gives
		\[
		-\frac4N\sum_{i,t}L_{i,t}^{\mathrm{corr}}(p^*)
		(\epsilon_t^k)^2\log(\epsilon_t^k)^{-1}+o(\Psi_k).
		\]
		The limiting relative scales and balance imply
		$(\epsilon_s^k/\epsilon_t^k)^{q_i}
		\to\rf_{p^*,s}^i(p_s^*)/\rf_{p^*,t}^i(p_t^*)$ for the minimal
		components. At $q_*=2$ the logarithms have ratio tending to one.
		Factoring the $t$-th scale proves \eqref{I:eq:sharp-global-preliminary}
		with precisely the coefficients in \eqref{I:def:tail-coefficients}.
		The balanced formula follows by factoring the mean scale.
	\end{proof}

	\subsection{Solutions at the critical parameter}
	
	\begin{corollary}\label{I:prop:stationary-leading-identity}
		Assume (H1)--(H2), (A1), and $h_i\in C^3(M)$ with $h_i>0$.
		Let $u^k$ be a sequence of blow-up solutions at
		$\rho^*\in\Gamma_N$. After a subsequence, its limiting configuration
		is a balanced critical pair $(p^*,r^*)$, and
		\begin{equation}\label{I:eq:stationary-leading-identity}
			0=-\mathcal L_t(p^*)\Psi_{\widehat m}(\widehat\lambda_t^k)
			+o(\Psi_k),\qquad t\in I'.
		\end{equation}
		With $\widehat L_k=N^{-1}\sum_t\log\widehat\lambda_t^k$,
		\begin{equation}\label{I:eq:stationary-balanced-leading-identity}
			0=-\mathcal L_\infty(p^*,r^*)
			\Psi_{\widehat m}(e^{\widehat L_k})
			+o\bigl(\Psi_{\widehat m}(e^{\widehat L_k})\bigr).
		\end{equation}
		Under (A2), no such sequence exists.
	\end{corollary}
	
	\begin{proof}
		Theorem~\ref{I:thm:path-free-coarse-matching} gives the balanced critical
		limit. Setting $\rho^k=\rho^*$ in
		\eqref{I:eq:sharp-global-preliminary} proves
		\eqref{I:eq:stationary-leading-identity}.
		Since $\log\widehat\lambda_t^k-\widehat L_k\to\delta_t^*$,
		\[
		\frac{\Psi_{\widehat m}(\widehat\lambda_t^k)}
		{\Psi_{\widehat m}(e^{\widehat L_k})}
		\longrightarrow e^{-(\widehat m-2)\delta_t^*}.
		\]
		Equation~\eqref{I:eq:L-compatible} gives
		\eqref{I:eq:stationary-balanced-leading-identity}. Division by the
		positive factor $\Psi_{\widehat m}(e^{\widehat L_k})$ contradicts (A2).
	\end{proof}
	
	\subsection{Refined estimates for the residual}
	
	\begin{proposition}
		Under the hypotheses of Theorem~\ref{I:thm:path-free-coarse-matching},
		use the parameters of Proposition~\ref{I:prop:sharp-bootstrap-closure}.
		Set
		\begin{equation*}
			\omega_k:=|\rho^k-\rho^*|+\Psi_k
			+(\widehat\lambda_1^k)^{-\gamma},
			\qquad 0<\gamma<\min\{1,\widehat m-2\}.
		\end{equation*}
		Choose $p\in(1,2)$ such that
		\begin{equation*}
			\theta_p:=2-\frac2p<\gamma.
		\end{equation*}
		For $\lambda_k:=\widehat\lambda_1^k$ and the residual in
		\eqref{I:eq:def-preliminary-residual},
		\begin{equation*}
			\max_i\|\widehat R_i^k\|_{L^p(M)}
			\le C\lambda_k^{\theta_p}\omega_k\longrightarrow0,
			\qquad
			\max_i\|\widehat R_i^k\|_{L^1(M)}\le C\omega_k.
		\end{equation*}
		Consequently, for $\zeta\in C^0(M)$,
		\[
		\left|\int_M\widehat R_i^k\zeta\,dV_g\right|
		\le C\omega_k\|\zeta\|_\infty.
		\]
	\end{proposition}
	
	\begin{proof}
		Lemma~\ref{I:lem:radial-family-derivatives} and the locally Lipschitz
		normalization give, for small fixed $\delta>0$,
		\[
		|e^{V_i(a;y)}-e^{v_i(y)}|
		\le C|a-\sigma|(1+\log(2+|y|))(1+|y|)^{-m_i+\delta}.
		\]
		By \eqref{I:eq:sharp-closed-masses}, $|a_t^k-\sigma|=O(\Psi_k)$.
		On $|y|\le C\lambda_k$, also
		\[
		(1+|y|)^{-(m_{i,t}^{V,k}-m_i)}=1+o(1),
		\qquad \Psi_k\log\lambda_k\to0.
		\]
		In a core, put $\widehat m_{k,t}:=\min_i m_{i,t}^{V,k}$. The actual source is
		$J_i^k(\epsilon y)e^{V_i^{k,t}(y)+\eta_i^k(y)}\,dy$, with
		$\epsilon=\widehat\lambda_t^{-1}$,
		$|J_i^k(\epsilon y)-1|\le C\epsilon|y|$, and, below four,
		\[
		|\eta_i^k(y)|\le C\epsilon|y|(1+|y|)^{2-\widehat m_{k,t}}
		+C\epsilon^2(1+|y|)^{4-\widehat m_{k,t}}\log(2+|y|).
		\]
		At four the corresponding bound is
		$C\epsilon|y|(1+|y|)^{-2}+C\epsilon^2\log^2(2+|y|)$.
		These errors tend uniformly to zero on the scaled fixed disk.
		The displacement contributes at most
		$C\epsilon(1+|y|)^{-m_i-1}$ to the density difference.
		
		The change of variables in the $L^p$ norm contributes
		$\lambda_k^{\theta_p}$. To estimate the remaining radial integrals, use
		\[
		\int_0^{C\lambda}(1+r)^{-a}r^b\log^d(2+r)\,r\,dr
		\le C_\delta\bigl(1+\lambda^{\max\{b+2-a,0\}+\delta}\bigr).
		\]
		The profile difference contributes $C\lambda_k^{\theta_p}\Psi_k$.
		The coefficient, displacement, and angular terms contribute at most
		$C\lambda_k^{\theta_p-\gamma}$ because
		$\gamma<\min\{1,\widehat m-2\}$.
		The refined radial term has boundary power
		$\lambda_k^{4-\widehat m-m_i}$ and satisfies the same bound;
		the strict choice of $\gamma$ absorbs the logarithms.
		
		On a fixed exterior, the circle and Green formulas give
		$s_i^k=O(\lambda_k^{-(m_i-2)})$. The projected tails and the mean
		correction contribute
		$O(|\rho^k-\rho^*|+\lambda_k^{-(\widehat m-2)})$.
		Adding these bounds proves the $L^p$ estimate. The same calculation
		with $p=1$ gives the $L^1$ bound. Finally,
		\eqref{I:eq:sharp-closed-masses} gives
		$|\rho^k-\rho^*|=O(\Psi_k)$, and $\theta_p<\gamma<\widehat m-2$
		implies $\lambda_k^{\theta_p}\omega_k\to0$.
	\end{proof}
	
	\subsection{Estimates for multiple blow-up points}
	\label{I:sec:multiple-coupling}
	
	For the preliminary parameters, set
	\[
	L_k:=\frac1N\sum_t\log\widehat\lambda_t^k,
	\qquad
	\nu_{i,t}^k:=\frac1{2\pi}
	\int_{U_{r_0}(\widehat p_t^k)}
	\rho_i^kh_ie^{\Theta_i^k}\,dV_g.
	\]
	We collect the moving-parameter estimates in the notation used below.
	
	\begin{proposition}
		Assume (H1)--(H2), (A1), and $h_i\in C^3(M)$ with $h_i>0$
		for every $i\in I$. Let $u^k$ be a sequence of blow-up solutions
		of \eqref{I:eq:solution}, where $\rho^k\to\rho^*\in\Gamma_N$
		is admissible in Definition~\ref{I:def:admissible-path}.
		For $N>1$, assume $\Lambda_{I,N}(\rho^k)>0$.
		After a subsequence and a permutation, the following hold.
		
		\begin{enumerate}[label=\textup{(\roman*)}]
			\item
			The concentration scales are comparable:
			\[
			C^{-1}\le\frac{\widehat\lambda_t^k}{\widehat\lambda_s^k}
			\le C,\qquad s,t\in I'.
			\]
			Moreover,
			\begin{equation}\label{I:eq:rho-rate}
				|\rho^k-\rho^*|\le C\Psi_k.
			\end{equation}
			
			\item
			The radial entire solutions $V^{k,t}$ in the local approximation,
			normalized by $\max_i V_i^{k,t}(0)=0$, satisfy
			\[
			\max_{i,t}\|V_i^{k,t}-v_i\|_{C^2(B_R)}\le C_R\Psi_k
			\qquad\text{for every fixed }R>0.
			\]
			The local masses and concentration points satisfy
			\[
			\max_{i,t}|\nu_{i,t}^k-\sigma_i|
			+|\nabla\ff(\widehat p^k)|\le C\Psi_k.
			\]
			
			\item
			The relative scales satisfy \eqref{I:eq:coefficient-balance}.
			
			\item
			For every $t\in I'$,
			\begin{equation}\label{I:eq:sharp-Lambda}
				\Lambda_{I,N}(\rho^k)
				=-\mathcal L_t(p^*)\Psi_{\widehat m}(\widehat\lambda_t^k)
				+o(\Psi_k).
			\end{equation}
		\end{enumerate}
	\end{proposition}
	
	\begin{proof}
		Scale comparability for $N>1$ follows from
		\cite[Theorem~1.1]{HuangZhang2022}.
		The source heights differ from \eqref{I:eq:height-scale} by bounded
		terms, so the concentration scales are comparable and the orders
		$O(\Psi_k)$ agree. The parameter bound is established in
		\cite[Lemma~3.2 and the proof of Theorem 1.2]{HuangZhang2022}.
		For $N=1$, \cite[Theorems~1.2--1.3]{LinZhang2013} gives
		$|\Lambda_{I,1}(\rho^k)|\le C\Psi_k$; transversality gives
		\eqref{I:eq:rho-rate}. These estimates also follow from
		\eqref{I:eq:sharp-closed-masses}.
		
		The local mass comparisons in \cite[Section~6]{LinZhang2013} and
		\cite[Theorem~1.2 and Section~3]{HuangZhang2022} give
		\[
		\nu_{i,t}^k=\frac{\rho_i^k}{2\pi N}+O(\Psi_k),
		\qquad
		\frac1{2\pi}\int_{\mathbb R^2}e^{V_i^{k,t}}\,dy
		=\nu_{i,t}^k+O(\Psi_k).
		\]
		Thus the entire masses differ from $\sigma$ by $O(\Psi_k)$.
		Lemma~\ref{I:lem:radial-family-derivatives} and the locally Lipschitz
		maximum-zero normalization give the $C^2(B_R)$ profile estimate.
		
		Proposition~\ref{I:prop:sharp-bootstrap-closure} gives the location
		estimate and the relative scale balance.
		Theorem~\ref{I:thm:sharp-global-mass-identity} gives
		\eqref{I:eq:sharp-Lambda}, with the coefficients defined in
		\eqref{I:def:tail-coefficients}.
	\end{proof}
	
	%------------------------------------------------------------------------------------------------------------------------------------------------
	\section{Orthogonal decomposition of blow-up solutions}
	\label{I:sec:entry-theorem}
	
	Throughout this section, $u^k$ satisfies the hypotheses of
	Theorem~\ref{I:thm:path-free-coarse-matching}. Let
	$(\widehat p^k,\widehat\lambda^k)$ be the preliminary parameters, and set
	\[
	\widehat\ell_t^k:=\log\widehat\lambda_t^k,
	\qquad \widehat L_k:=\frac1N\sum_t\widehat\ell_t^k.
	\]
	Fix $1<p<1+\beta/2$ and
	$\kappa_p=(\beta-2(p-1))/p>0$. Proposition~\ref{I:prop:coarse-source-entry}
	and Sobolev embedding give
	\[
	\|u^k-U_{\widehat p^k,\widehat\lambda^k,0}\|_{\mathcal H}
	+\|u^k-U_{\widehat p^k,\widehat\lambda^k,0}\|_\infty
	\le C\|u^k-U_{\widehat p^k,\widehat\lambda^k,0}\|_{W^{2,p}(M)^n}
	\le Ce^{-\kappa_p\widehat L_k}.
	\]
	
	\subsection{Minimizing parameters}
	
	\begin{lemma}\label{I:lem:modulation-estimate}
		For all sufficiently large $k$, problem
		\eqref{I:eq:exact-parameter-minimization}, with
		$\bar z=(\widehat p^k,\widehat\ell^k,0)$, has a unique ordered minimizer
		$(p^k,\ell^k,b^k)$. Set
		\[
		\lambda^k:=\lambda(\ell^k),\qquad
		w^k:=u^k-U_{p^k,\lambda^k,b^k},\qquad
		L_k:=\frac1N\sum_t\ell_t^k.
		\]
		Then $w^k\in E^N_{p^k,\lambda^k,b^k}$ and
		\begin{equation}\label{I:eq:modulation-rate}
			\begin{aligned}
				&\max_t\widehat\lambda_t^k d_g(p_t^k,\widehat p_t^k)
				+\max_t|\ell_t^k-\widehat\ell_t^k|
				+\sqrt{\widehat L_k}|b^k|+\|w^k\|_{\mathcal H}\\
				&\hspace{15mm}\le
				C\|u^k-U_{\widehat p^k,\widehat\lambda^k,0}\|_{\mathcal H}
				\le Ce^{-\kappa_p\widehat L_k}.
			\end{aligned}
		\end{equation}
		The minimizing parameters are independent of the preliminary
		representation after relabeling. Moreover,
		\begin{equation}\label{I:eq:b-rate-intermediate}
			|b^k|\le C\widehat L_k^{-1/2}e^{-\kappa_p\widehat L_k},
			\qquad
			L_k|b^k|\le C\sqrt{\widehat L_k}e^{-\kappa_p\widehat L_k}\longrightarrow0.
		\end{equation}
	\end{lemma}
	
	\begin{proof}
		We first verify the gradient condition in
		$V^\eta_{\rho^*}(N,\eps)$. By \eqref{I:eq:def-preliminary-residual},
		\[
		\rho_i^k(h_i e^{\Theta_i^k}-1)
		=-\Delta_g\sum_t PW_{i,\widehat p_t^k,\widehat\lambda_t^k}
		+\widehat R_i^k.
		\]
		Proposition~\ref{I:prop:coarse-source-entry} and
		\eqref{I:eq:gram-bb} therefore imply
		\[
		\|h_i e^{\Theta_i^k}-1\|_{(\oH)^*}
		\le C\left(\left\|\sum_t
		PW_{i,\widehat p_t^k,\widehat\lambda_t^k}\right\|_{\oH}
		+\|\widehat R_i^k\|_{L^p(M)}\right)
		\le C(1+\sqrt{\widehat L_k}).
		\]
		Since $DJ_{\rho^k}(u^k)=0$, for $v\in\mathcal H$,
		\[
		DJ_{\rho^*}(u^k)[v]
		=\sum_i(\rho_i^k-\rho_i^*)
		\int_M(h_i e^{\Theta_i^k}-1)v_i\,dV_g.
		\]
		Hence \eqref{I:eq:path-free-mean-comparability} gives
		\[
		\|DJ_{\rho^*}(u^k)\|_{\mathcal H^*}
		\le C|\rho^k-\rho^*|(1+\sqrt{\widehat L_k})\longrightarrow0.
		\]
		Together with the preliminary approximation, separation, and scale
		comparability, this proves that $u^k\in V^\eta_{\rho^*}(N,\eps)$ for
		every fixed small $\eps>0$ and all large $k$.
		Lemma~\ref{I:lem:unique} gives the unique interior minimizer, its
		independence of the preliminary representation, and
		$w^k\in E^N_{p^k,\lambda^k,b^k}$.
		
		By minimality,
		\[
		\|w^k\|_{\mathcal H}
		\le\|u^k-U_{\widehat p^k,\widehat\lambda^k,0}\|_{\mathcal H},
		\qquad
		\|U_{p^k,\lambda^k,b^k}-U_{\widehat p^k,\widehat\lambda^k,0}\|_{\mathcal H}
		\le2\|u^k-U_{\widehat p^k,\widehat\lambda^k,0}\|_{\mathcal H}.
		\]
		Apply the lower bound in \eqref{I:eq:modulation-bilipschitz} in the
		normalized coordinates of Lemma~\ref{I:lem:unique}. This proves
		\eqref{I:eq:modulation-rate}. In particular,
		$|L_k-\widehat L_k|\le Ce^{-\kappa_p\widehat L_k}$, which gives
		\eqref{I:eq:b-rate-intermediate}.
	\end{proof}
	
	\subsection{Local mass balance}
	
	For $i\in I$ and $t\in I'$, define
	\begin{equation}\label{I:def:full-local-mass}
		\mathfrak M_{i,t}(p,\lambda,b):=
		2\pi\sigma_{i,t}^{\alpha}
		\exp\!\left(-\sum_j\alpha^{i,j,t}\fJ_j\right)
		(\rf_p^i\rg_{p,\alpha}^i)(p_t)
		\lambda_t^{q_i+\sigma_i b_{i,t}}.
	\end{equation}
	For a minimizing decomposition $u=U_{p,\lambda,b}+w$, set
	\begin{equation}\label{I:def:tau_ijt}
		\tau_{i,t}:=1-
		\frac{N\mathfrak M_{i,t}(p,\lambda,b)}{\int_M h_i e^{u_i}\,dV_g},
		\qquad
		\widehat\tau_{i,t}:=1-
		\frac{N\mathfrak M_{i,t}(p,\lambda,b)}
		{\sum_s\mathfrak M_{i,s}(p,\lambda,b)}.
	\end{equation}
	Thus $\sum_t\widehat\tau_{i,t}=0$. For the sequence $u^k$, write
	$\mathfrak M_{i,t}^k:=\mathfrak M_{i,t}(p^k,\lambda^k,b^k)$,
	$\tau_{i,t}^k$, and $\widehat\tau_{i,t}^k$.
	
	\begin{lemma}\label{I:lem:centered-entry}
		The minimizing decomposition satisfies
		\[
		\max_i\left|
		\frac{\int_M h_i e^{u_i^k}\,dV_g}{\sum_t\mathfrak M_{i,t}^k}-1
		\right|
		+\max_{i,t}\bigl(|\tau_{i,t}^k|+|\widehat\tau_{i,t}^k|\bigr)
		\longrightarrow0.
		\]
	\end{lemma}
	
	\begin{proof}
		On a fixed disk about $\widehat p_t^k$ where $\chi_{\widehat p_t^k}=1$,
		the projected expansion
		\eqref{I:eq:exact-projected-expansion} gives
		\[
		h_i e^{U_{i,\widehat p^k,\widehat\lambda^k,0}}
		=(\widehat\lambda_t^k)^{m_i}
		e^{v_i(\widehat\lambda_t^k y_{\widehat p_t^k})-\fJ_i}
		\rf_{\widehat p^k,t}^i(x)(1+o(1)).
		\]
		After rescaling, the integrand divided by its smooth coefficient is
		bounded by $Ce^{v_i}$, whose integral is $2\pi\sigma_i$.
		The complement of the disks contributes $O(1)$, whereas $q_i>0$.
		Dominated convergence and the preliminary $C^0$ approximation yield
		\[
		\int_M h_i e^{u_i^k}\,dV_g
		=2\pi\sigma_i e^{-\fJ_i}
		\left(\sum_t\rf_{\widehat p^k,t}^i(\widehat p_t^k)
		(\widehat\lambda_t^k)^{q_i}\right)(1+o(1)).
		\]
		
		For small $b$, the logarithmic moments defining
		$D_b\sigma_{i,t}^{\alpha}$ are uniformly integrable. Differentiating
		\eqref{I:def:full-local-mass} on the parameter box gives
		\[
		|D_b\log\mathfrak M_{i,t}|\le C(1+\widehat L_k),
		\qquad
		|D_p\log\mathfrak M_{i,t}|+|D_\ell\log\mathfrak M_{i,t}|\le C.
		\]
		Since
		$\mathfrak M_{i,t}(p,\lambda,0)
		=2\pi\sigma_i e^{-\fJ_i}\rf_{p,t}^i(p_t)\lambda_t^{q_i}$,
		\eqref{I:eq:modulation-rate} implies
		\[
		\max_{i,t}\left|\log
		\frac{\mathfrak M_{i,t}^k}
		{2\pi\sigma_i e^{-\fJ_i}
			\rf_{\widehat p^k,t}^i(\widehat p_t^k)(\widehat\lambda_t^k)^{q_i}}
		\right|
		\le C\sqrt{\widehat L_k}e^{-\kappa_p\widehat L_k}\longrightarrow0.
		\]
		These two formulas prove the total-mass assertion.
		
		By Theorem~\ref{I:thm:path-free-coarse-matching},
		$\widehat p^k\to p^*$ and
		$\widehat\ell_t^k-\widehat L_k\to\delta_t^*$, with
		$\rf_{p^*,t}^i(p_t^*)e^{q_i\delta_t^*}$ independent of $t$.
		It follows that
		\[
		\frac{\mathfrak M_{i,t}^k}{\sum_s\mathfrak M_{i,s}^k}
		=\frac{\rf_{\widehat p^k,t}^i(\widehat p_t^k)
			e^{q_i(\widehat\ell_t^k-\widehat L_k)}}
		{\sum_s\rf_{\widehat p^k,s}^i(\widehat p_s^k)
			e^{q_i(\widehat\ell_s^k-\widehat L_k)}}+o(1)
		\longrightarrow\frac1N.
		\]
		Thus $\widehat\tau_{i,t}^k\to0$. Finally,
		\[
		\tau_{i,t}^k
		=1-\frac{\sum_s\mathfrak M_{i,s}^k}
		{\int_M h_i e^{u_i^k}\,dV_g}(1-\widehat\tau_{i,t}^k)
		\longrightarrow0.
		\]
	\end{proof}
	
	\subsection{Asymptotic decomposition}
	
	\begin{theorem}
		Under the hypotheses of Theorem~\ref{I:thm:path-free-coarse-matching},
		after a subsequence and relabeling, the solutions have a unique minimizing
		decomposition, up to permutation,
		\begin{equation}\label{I:eq:entry-decomposition}
			u^k=U_{p^k,\lambda^k,b^k}+w^k,
			\qquad w^k\in E^N_{p^k,\lambda^k,b^k}.
		\end{equation}
		With $L_k=N^{-1}\sum_t\log\lambda_t^k$, one has
		\[
		\min_t\lambda_t^k\longrightarrow\infty,
		\qquad C^{-1}\le\frac{\lambda_t^k}{\lambda_s^k}\le C,
		\qquad
		L_k|b^k|+\|w^k\|_{\mathcal H}
		+\max_{i,t}\bigl(|\tau_{i,t}^k|+|\widehat\tau_{i,t}^k|\bigr)
		\longrightarrow0.
		\]
		Moreover,
		\[
		p^k\longrightarrow p^*,\qquad
		\log\lambda_t^k-L_k\longrightarrow\delta_t^*,\qquad
		\sum_t\delta_t^*=0,
		\]
		where
		\[
		\nabla\mathcal F_{\rho^*,h}^N(p^*)=0,
		\qquad
		\rf_{p^*,t}^i(p_t^*)e^{q_i\delta_t^*}
		=\rf_{p^*,s}^i(p_s^*)e^{q_i\delta_s^*}
		\quad(i\in I,\ s,t\in I').
		\]
		Writing $\delta_t^*=\sum_{\alpha=1}^{N-1}e_t^\alpha r_\alpha^*$,
		the parameter defect satisfies, for every $t$,
		\[
		\begin{aligned}
			\Lambda_{I,N}(\rho^k)
			&=-\mathcal L_t(p^*)\Psi_{\widehat m}(\lambda_t^k)
			+o\bigl(\Psi_{\widehat m}(\lambda_t^k)\bigr)\\
			&=-\mathcal L_\infty(p^*,r^*)\Psi_{\widehat m}(e^{L_k})
			+o\bigl(\Psi_{\widehat m}(e^{L_k})\bigr).
		\end{aligned}
		\]
	\end{theorem}
	
	\begin{proof}
		Lemmas~\ref{I:lem:modulation-estimate} and \ref{I:lem:centered-entry}
		give the decomposition and smallness estimates. By
		\eqref{I:eq:modulation-rate},
		\[
		d_g(p_t^k,\widehat p_t^k)\longrightarrow0,
		\qquad \log\frac{\lambda_t^k}{\widehat\lambda_t^k}\longrightarrow0,
		\qquad L_k-\widehat L_k\longrightarrow0.
		\]
		The limits and balance conditions therefore follow from
		Theorem~\ref{I:thm:path-free-coarse-matching}.
		In both decay regimes,
		\[
		\frac{\Psi_{\widehat m}(\lambda_t^k)}
		{\Psi_{\widehat m}(\widehat\lambda_t^k)}\longrightarrow1,
		\qquad
		\frac{\Psi_{\widehat m}(\lambda_t^k)}
		{\Psi_{\widehat m}(e^{L_k})}
		\longrightarrow e^{-(\widehat m-2)\delta_t^*}.
		\]
		The first limit transfers the identity of
		Theorem~\ref{I:thm:sharp-global-mass-identity} to the minimizing scales.
		The second and Proposition~\ref{I:prop:balanced-coefficient-compatibility}
		give its balanced form.
	\end{proof}

	%------------------------------------------------------------------------
	\section{Uniform capture, separation, and finiteness}
	\label{I:sec:uniform-capture}
	
	We first prove finiteness of the balanced critical configurations, then
	choose disjoint neighborhoods that contain all sufficiently large
	solutions along a fixed admissible path.
	
	\subsection{Finiteness of balanced critical configurations}
	
	\begin{proposition}\label{I:prop:finiteness-balanced-orbits}
		For $N\ge2$, there exists $d_*>0$ such that every critical point of
		$\ff$ satisfies
		\[
		\min_{s\ne t}d_g(p_s,p_t)\ge d_*.
		\]
		Under \textup{(A1)}, the ordered set $\mathcal K_{\mathrm{bal}}$ in
		\eqref{I:def:balanced-critical-set} is finite.
	\end{proposition}
	
	\begin{proof}
		Suppose that critical configurations $p^k$ satisfy
		$d_k:=\min_{s\ne t}d_g(p_s^k,p_t^k)\to0$.
		After a subsequence, all pairwise distance ratios divided by $d_k$
		have limits in $[1,\infty]$. The labels whose mutual ratios are bounded
		form clusters. Choose the cluster $J$ containing a closest pair. Then
		\[
		|J|\ge2,\qquad
		\max_{s,t\in J}d_g(p_s^k,p_t^k)=O(d_k),\qquad
		\frac{d_g(p_s^k,p_t^k)}{d_k}\longrightarrow\infty
		\quad(s\in J,\ t\notin J).
		\]
		Use a fixed isothermal chart about the limiting cluster point, write
		$x_t^k$ for the coordinates of $p_t^k$ with $t\in J$, and set
		\[
		\bar x_k:=\frac1{|J|}\sum_{t\in J}x_t^k,
		\qquad
		X_t^k:=
		\begin{cases}
			x_t^k-\bar x_k,&t\in J,\\
			0,&t\notin J.
		\end{cases}
		\]
		By \eqref{I:def:reduce-F}, the singular interaction within $J$ is
		\[
		-\left(\sum_i\rho_i^*m_i\right)
		\sum_{\substack{s<t\\s,t\in J}}\log|x_s-x_t|.
		\]
		The smooth terms contribute $O(d_k)$ to $D\ff(p^k)[X^k]$;
		interactions with labels outside $J$ contribute
		\[
		O\!\left(\sum_{\substack{s\in J\\t\notin J}}
		\frac{d_k}{d_g(p_s^k,p_t^k)}\right)=o(1).
		\]
		For each internal pair,
		\[
		\frac{x_s^k-x_t^k}{|x_s^k-x_t^k|^2}\cdot X_s^k
		+\frac{x_t^k-x_s^k}{|x_s^k-x_t^k|^2}\cdot X_t^k=1.
		\]
		Consequently,
		\[
		0=D\ff(p^k)[X^k]
		=-\left(\sum_i\rho_i^*m_i\right)\binom{|J|}{2}+o(1),
		\]
		which is impossible. Thus the critical set is compact in
		$M^N\setminus\Delta$. For $N=1$, compactness follows directly from
		compactness of $M$.
		
		Under (A1), this critical set is discrete and hence finite.
		For each fixed $p$, the balance equations
		\[
		q_i(\delta_t-\delta_s)
		=-\log\frac{\rf_{p,t}^i(p_t)}{\rf_{p,s}^i(p_s)},
		\qquad \sum_t\delta_t=0,
		\]
		have at most one solution because $q_i>0$.
		Thus $\mathcal K_{\mathrm{bal}}$ is finite.
	\end{proof}
	
	\subsection{Bubbling neighborhoods}
	
	Assume (A1) and first suppose $\mathcal K_{\mathrm{bal}}\ne\varnothing$.
	Choose $\eta>0$ so small that the closed ordered center neighborhoods
	\[
	\left\{p:\max_t d_g(p_t,p_t^*)\le2\eta\right\},
	\qquad (p^*,r^*)\in\mathcal K_{\mathrm{bal}},
	\]
	are pairwise disjoint and lie inside the configuration neighborhoods used
	in Lemma~\ref{I:lem:unique}. Such a choice is possible by
	Proposition~\ref{I:prop:finiteness-balanced-orbits} and uniqueness of
	the balance vector at each location.
	
	For $(p^*,r^*)\in\mathcal K_{\mathrm{bal}}$, let
	$\delta^*$ be given by \eqref{I:def:delta-star}.
	Define $\mathcal U_{\eta,L_0}([p^*,r^*])$ as the set of functions
	whose minimizing decomposition \eqref{I:eq:entry-decomposition}, written as
	\[
	u=U_{p,\lambda,b}+w,\qquad w\in E^N_{p,\lambda,b},
	\qquad L:=\frac1N\sum_t\log\lambda_t,
	\]
	admits an ordering satisfying
	\begin{equation}\label{I:eq:capture-neighborhood}
		\begin{aligned}
			&L>L_0,\qquad
			\max_t d_g(p_t,p_t^*)<\eta,\qquad
			\max_t|\log\lambda_t-L-\delta_t^*|<\eta,\\
			&L|b|+\|w\|_{\mathcal H}
			+\max_{i,t}|\widehat\tau_{i,t}|<\eta.
		\end{aligned}
	\end{equation}
	Relabeling gives
	\[
	U_{\pi p,\pi\lambda,\pi b}=U_{p,\lambda,b},
	\qquad
	(\pi r)_\alpha=\sum_t e_t^\alpha\delta_{\pi^{-1}(t)}.
	\]
	Hence this set depends only on the orbit $[p^*,r^*]$.
	
	\begin{lemma}\label{I:prop:uniform-separation}
		For sufficiently small $\eta$ and sufficiently large $L_0$, these sets
		are open in $\mathcal H$, and their closures in $\mathcal H$ are pairwise
		disjoint for distinct balanced orbits.
	\end{lemma}
	
	\begin{proof}
		The normalized Hessian in the proof of Lemma~\ref{I:lem:unique} satisfies
		\[
		D_\xi^2\!\left(\frac12\|u-U(\xi)\|_{\mathcal H}^2\right)[h,h]
		\ge(c-C\|w\|_{\mathcal H})|h|^2
		\ge\frac c2|h|^2
		\]
		at the minimizing parameters. The same least squares argument applies
		uniformly on the present neighborhoods, since
		\[
		\|u-U_{p,\lambda,0}\|_{\mathcal H}
		\le\|w\|_{\mathcal H}+C\sqrt L\,|b|
		\le\eta\left(1+\frac C{\sqrt{L_0}}\right).
		\]
		Lemma~\ref{I:lem:preliminary-identification} puts overlapping
		representations in one small normalized chart after relabeling.
		The positive Hessian gives uniqueness there, and the implicit function
		theorem gives a locally $C^1$ map $u\mapsto(p,\lambda,b,w)$.
		The strict inequalities in \eqref{I:eq:capture-neighborhood} define
		open sets.
		
		To pass to closures, let $u^k\to u$ in $\mathcal H$ with $u^k$ in one
		such neighborhood. On these parameter sets, the energy identity
		\eqref{I:eq:bubble-energy} gives
		\[
		\|U_{p,\lambda,b}\|_{\mathcal H}^2
		=2\pi\sum_{i,t}(m_i+\sigma_i b_{i,t})^2\log\lambda_t+O(1)
		\ge cL-C.
		\]
		Since $\|U_{p^k,\lambda^k,b^k}\|_{\mathcal H}
		\le\|u^k\|_{\mathcal H}+\eta$, the means $L_k$ are bounded.
		The centers, relative scales, and amplitudes then range in a compact
		parameter set. After a subsequence, the parameters converge and the
		orthogonality relations pass to the limit. The positive Hessian
		identifies the limit as the unique ordered minimizer in the larger
		chart. It satisfies \eqref{I:eq:capture-neighborhood} with non-strict
		inequalities.
		
		If $u$ belonged to the closures of two neighborhoods, align their
		limiting minimizing parameters. Their common center vector would satisfy
		\[
		\max_t d_g(p_t^*,p_t^{\prime *})
		\le\max_t d_g(p_t^*,p_t)+\max_t d_g(p_t,p_t^{\prime *})
		\le2\eta
		\]
		after relabeling one representative. The choice of the ordered center
		neighborhoods forces $p^*=p^{\prime *}$, and uniqueness of scale balance
		then gives the same balanced orbit.
	\end{proof}
	
	\subsection{Uniform capture and compactness}
	
	\begin{theorem}
		Assume (A1). Fix a $C^1$ path $\rho(s)$ with $\rho(0)=\rho^*$ whose
		derivative at zero has nonzero components of one common sign.
		There exist $\eta>0$, $L_0>1$, $s_0>0$, and $R_0>0$ such that every
		solution at parameter $\rho(s)$ satisfying
		\[
		0<|s|<s_0,\qquad
		\Lambda_{I,N}(\rho(s))>0\quad\text{if }N>1,
		\qquad \max_i\|u_i\|_{L^\infty(M)}>R_0
		\]
		belongs to exactly one of the finitely many neighborhoods
		\[
		\mathcal U_{\eta,L_0}([p^*,r^*]),
		\qquad [p^*,r^*]\in\mathcal K_{\mathrm{bal}}/\mathfrak S_N.
		\]
		Their closures are pairwise disjoint, and the minimizing decomposition
		is unique up to permutation.
	\end{theorem}
	
	\begin{proof}
		Suppose first that $\mathcal K_{\mathrm{bal}}\ne\varnothing$, and choose
		$\eta,L_0$ as in Lemma~\ref{I:prop:uniform-separation}.
		If no $s_0,R_0$ worked, there would be solutions $u^k$ at $\rho(s_k)$
		outside all these neighborhoods, with
		\[
		s_k\to0,\qquad
		\max_i\|u_i^k\|_{L^\infty(M)}\to\infty.
		\]
		The sequence $\rho(s_k)$ is admissible. After a subsequence,
		Theorem~\ref{I:thm:intro-coarse-entry} and
		Lemma~\ref{I:lem:centered-entry} give a balanced pair $(p^*,r^*)$ and
		\[
		L_k\to\infty,\qquad p^k\to p^*,\qquad
		\log\lambda_t^k-L_k\to\delta_t^*,\qquad
		L_k|b^k|+\|w^k\|_{\mathcal H}
		+\max_{i,t}|\widehat\tau_{i,t}^k|\to0.
		\]
		Thus \eqref{I:eq:capture-neighborhood} holds for large $k$, a
		contradiction. If $\mathcal K_{\mathrm{bal}}$ is empty, the same
		sequence argument contradicts the existence of its balanced limit and
		gives a uniform bound for all solutions on the stated part of the path.
		Uniqueness of the neighborhood follows from
		Lemma~\ref{I:prop:uniform-separation}.
	\end{proof}
	
	\begin{corollary}
		Assume (A1)--(A2). The solution set at $\rho^*$ is bounded in
		$L^\infty(M)^n$ and compact in $C^{2,\alpha}(M)^n$ for every
		$0<\alpha<1$.
	\end{corollary}
	
	\begin{proof}
		An unbounded sequence of solutions at $\rho^*$ would contradict
		Corollary~\ref{I:prop:stationary-leading-identity} under (A2).
		Elliptic regularity gives precompactness in $C^{2,\alpha}(M)^n$;
		passing to the limit in the equation gives closedness.
	\end{proof}
	
	%---------------------------------------------------------------------

	%---------------------------
	\appendix
	\section{Projected bubbles and parameter estimates}
	
	We prove the estimates used in the least squares decomposition of
	Section~\ref{I:sec:bubbling-chart} and the closure argument of
	Section~\ref{I:sec:uniform-capture}. Fix the entire profile $v$ and the
	normalized chart families. Centers range in a compact subset of
	$M^N\setminus\Delta$, with a common cutoff radius $0<r_0<1$ and comparable
	scales. Write $q_*:=\widehat m-2\in(0,2]$ and $\ell:=\log\lambda$.
	
	\begin{lemma}\label{I:lem:ex-bubble}
		For one center $p$ and scale $\lambda\ge2/r_0$, set
		\[
		c_{p,i}:=-\int_M\chi_p
		\bigl(v_i(\lambda y_p)+m_i\log\lambda-\fJ_i
		+m_i\log|y_p|\bigr)\,dV_g.
		\]
		Then
		\begin{equation}\label{I:eq:exact-projected-expansion}
			PW^*_{i,p,\lambda}
			=\chi_p\bigl(v_i(\lambda y_p)+m_i\log\lambda-\fJ_i\bigr)
			+2\pi m_iH^g(\cdot,p)+c_{p,i}+R_{p,i}.
		\end{equation}
		Here $\int_MR_{p,i}\,dV_g=0$. For $|a|+j\le2$, uniformly in $p$ and
		$\lambda$,
		\begin{equation}\label{I:eq:projected-remainder-derivatives}
			\|D_p^a\partial_\ell^jR_{p,i}\|_{C^1(M)}
			\le C\lambda^{-q_*},
			\qquad
			|D_p^a\partial_\ell^jc_{p,i}|
			\le C\Psi_{\widehat m}(\lambda).
		\end{equation}
		Here $D_p$ denotes derivatives in a fixed center chart. The unstarred
		projections satisfy $PW_i=\sum_j a^{ij}PW_j^*$.
	\end{lemma}
	
	\begin{proof}
		Set $w_i(r):=v_i(r)+m_i\log r-\fJ_i$.
		The radial identity \eqref{I:eq:radial-volterra} and the equations
		\[
		rw_i'(r)=\sum_j a_{ij}\int_r^\infty s e^{v_j(s)}\,ds,
		\qquad
		(r\partial_r)^2w_i(r)=-r^2\sum_j a_{ij}e^{v_j(r)}
		\]
		give
		\[
		|(r\partial_r)^j w_i(r)|\le C
		\begin{cases}
			r^{-q_*},&r\ge1,\\
			1+|\log r|,&0<r\le1,
		\end{cases}
		\quad 0\le j\le3.
		\]
		Using $H^g=G^g+(2\pi)^{-1}\chi_p\log|y_p|$, the first two terms of
		\eqref{I:eq:exact-projected-expansion} equal
		$2\pi m_iG^g(\cdot,p)+\chi_pw_i(\lambda|y_p|)$.
		Their mean is $-c_{p,i}$, so $\int_MR_{p,i}=0$.
		The defining Poisson equations give
		\[
		\begin{aligned}
			-\Delta_gR_{p,i}={}&
			2\langle\nabla_g\chi_p,\nabla_g w_i(\lambda|y_p|)\rangle_g
			+(\Delta_g\chi_p)w_i(\lambda|y_p|)\\
			&+2\pi m_i-\sum_j a_{ij}
			\int_{\mathbb R^2}\chi\!\left(\frac{|y|}{\lambda r_0}\right)
			e^{v_j(y)}\,dy.
		\end{aligned}
		\]
		The cutoff derivatives are supported on a fixed annulus, and
		$2\pi m_i=\sum_j a_{ij}\int_{\mathbb R^2}e^{v_j}$.
		The radial bounds therefore yield
		\[
		\|D_p^a\partial_\ell^j\Delta_gR_{p,i}\|_\infty
		\le C\lambda^{-q_*},\qquad |a|+j\le2.
		\]
		The mean-zero elliptic estimate proves the first bound in
		\eqref{I:eq:projected-remainder-derivatives}.
		
		For the scalar term, change variables to the fixed chart disk:
		\[
		c_{p,i}=-\int_{B_{r_0}}
		\chi(|y|/r_0)w_i(\lambda|y|)e^{\widehat\varphi_p(y)}\,dy.
		\]
		Center derivatives act only on the smooth Jacobian. The same radial
		bounds give, for $|a|+j\le2$,
		\[
		|D_p^a\partial_\ell^jc_{p,i}|
		\le C\left[
		\int_0^{\lambda^{-1}}(1+|\log(\lambda r)|)r\,dr
		+\lambda^{-q_*}\int_{\lambda^{-1}}^{r_0}r^{1-q_*}\,dr
		\right]
		\le C\Psi_{\widehat m}(\lambda).
		\]
		This proves the remaining bound.
	\end{proof}
	
	\begin{lemma}\label{I:lem:uniform-C2-bubble-map}
		For $\bar z=(\bar p,\bar\ell,0)$, put
		$\bar\lambda=\lambda(\bar\ell)$ and
		$\bar L=N^{-1}\sum_t\bar\ell_t$. Use the normalized coordinates
		\[
		\begin{aligned}
			\xi_{\bar z}(p,\ell,b)
			&:={\left(
				(\bar\lambda_t y_{\bar p_t}(p_t))_t,
				(\ell_t-\bar\ell_t)_t,
				(\sqrt{\bar L}\,b_{i,t})_{i,t}\right)},\\
			\Theta_{\bar z}(\xi)&:=U_{p,\lambda(\ell),b}.
		\end{aligned}
		\]
		For each fixed small $r>0$, there are $L_0,C_r$ such that
		\[
		\|D_\xi\Theta_{\bar z}\|_{\mathrm{op}}
		+\|D_\xi^2\Theta_{\bar z}\|_{\mathrm{op}}
		\le C_r
		\qquad(\bar L\ge L_0,\ |\xi|\le r),
		\]
		as operators into $\mathcal H$, uniformly in $\bar z$.
	\end{lemma}
	
	\begin{proof}
		The radial equation and the preceding bounds imply, for $|y|\ge1$,
		\[
		|D^av_i(y)|\le C|y|^{-|a|}\quad(1\le|a|\le3),
		\qquad
		|\nabla(y\cdot\nabla v_i+m_i)|\le C|y|^{-1-q_*}.
		\]
		Rescaling in each core and using Lemma~\ref{I:lem:ex-bubble} gives
		\[
		\|PW^*_{i,p,\lambda}\|_{\oH}\le C\sqrt{\log\lambda},
		\qquad
		\lambda^{-|a|}\|D_p^a\partial_\ell^jPW^*_{i,p,\lambda}\|_{\oH}
		\le C\quad(1\le|a|+j\le2).
		\]
		Indeed the gradients of the translation and dilation derivatives are
		square integrable on $\mathbb R^2$. Derivatives of the smooth chart
		coefficients contribute at most $C\lambda^{-1}\sqrt{\log\lambda}$.
		On the cutoff annulus, differentiate the combined expression
		\[
		\chi_p\bigl(v_i(\lambda y_p)+m_i\ell-\fJ_i\bigr)
		=-m_i\chi_p\log|y_p|+\chi_pw_i(\lambda|y_p|).
		\]
		Its normalized derivatives are bounded; the regular Green terms and
		$R_{p,i}$ satisfy the same bounds. The scalars $c_{p,i}$ have zero
		Dirichlet contribution.
		
		Since $PW_i=\sum_j a^{ij}PW_j^*$, these estimates also hold for $PW_i$.
		The map $U$ is affine in $b$, and
		\[
		\bar L^{-1/2}\|PW_{i,p_t,\lambda_t}\|_{\oH}\le C,
		\qquad D_b^2U=0.
		\]
		The chain rule, with $\lambda_t/\bar\lambda_t$ bounded above and below,
		proves the assertion, including the mixed amplitude derivatives.
	\end{proof}
	
	\begin{lemma}
		The amplitude directions satisfy
		\begin{equation}\label{I:eq:gram-bb}
			\langle X^b_{i,t},X^b_{j,s}\rangle_{\mathcal H}
			=\delta_{ij}\delta_{st}\,2\pi\sigma_i^2\log\lambda_t+O(1).
		\end{equation}
		Moreover, uniformly for bounded $b$,
		\begin{equation}\label{I:eq:bubble-energy}
			\|U_{p,\lambda,b}\|_{\mathcal H}^2
			=2\pi\sum_{i,t}(m_i+\sigma_i b_{i,t})^2\log\lambda_t+O(1).
		\end{equation}
	\end{lemma}
	
	\begin{proof}
		Set $\psi_i:=\sum_j a^{ij}v_j$. Since $A^{-1}m=\sigma$,
		\[
		\nabla\psi_i(y)=-\sigma_i\frac{y}{|y|^2}
		+O(|y|^{-1-q_*})\quad(|y|\ge1).
		\]
		The projected expansion gives
		\[
		\langle PW_{i,p_t,\lambda_t},PW_{j,p_s,\lambda_s}\rangle_{\oH}
		=\delta_{st}\,2\pi\sigma_i\sigma_j\log\lambda_t+O(1).
		\]
		For $s=t$, the logarithmic annulus gives the leading term, and the
		core and error terms are bounded. For $s\ne t$, separation bounds the
		projection on the support of the other bubble's Poisson density.
		Since distinct components are orthogonal in $\mathcal H$, these scalar
		products imply \eqref{I:eq:gram-bb}. Expanding the norm of
		$U_i=\sum_{j,t}(a_{ij}+\delta_{ij}b_{i,t})PW_{j,p_t,\lambda_t}$
		and using $\sum_j a_{ij}\sigma_j=m_i$ gives
		\eqref{I:eq:bubble-energy}.
	\end{proof}
	
	\begin{lemma}\label{I:lem:tangent-gram}
		Let $\Theta_{\bar z}$ and $\xi$ be as in
		Lemma~\ref{I:lem:uniform-C2-bubble-map}. There exist $r,c,C,L_0>0$ such
		that
		\[
		c|h|^2\le\|D_\xi\Theta_{\bar z}(\xi)h\|_{\mathcal H}^2
		\le C|h|^2
		\qquad(\bar L\ge L_0,\ |\xi|\le r)
		\]
		for every parameter vector $h$, uniformly in $\bar z$.
	\end{lemma}
	
	\begin{proof}
		At $\xi=0$, one has $(p,\ell,b)=(\bar p,\bar\ell,0)$.
		The columns of $D_\xi\Theta_{\bar z}(0)$ are
		\[
		\bar\lambda_t^{-1}X^p_{t,l},\qquad
		\bar\lambda_tX_t^\lambda,\qquad
		\bar L^{-1/2}X^b_{i,t},
		\]
		evaluated at $\bar z$, with position derivatives in the
		$y_{\bar p_t}$ coordinates.
		The position and scale columns have rescaled profiles
		$-\partial_l v_i$ and $y\cdot\nabla v_i+m_i$, whose gradients satisfy
		\[
		|\nabla\partial_l v_i(y)|\le C|y|^{-2},\qquad
		|\nabla(y\cdot\nabla v_i+m_i)|\le C|y|^{-1-q_*}
		\quad(|y|\ge1).
		\]
		Their Dirichlet norms are finite and positive: the dilation profile
		equals $m_i$ at zero and tends to zero at infinity, and
		$v_i''(0)=-\tfrac12\sum_j a_{ij}e^{v_j(0)}<0$.
		The two translation fields and the dilation field are mutually
		orthogonal by radial symmetry; separation makes the different-bubble
		products tend to zero.
		Their products with unnormalized amplitude directions are $O(1)$,
		by the integrable radial tails $r^{-2}$ and $r^{-1-q_*}$.
		After amplitude normalization these mixed products tend to zero, while
		\eqref{I:eq:gram-bb} and $\log\bar\lambda_t/\bar L\to1$ give positive
		limiting amplitude norms. Thus the asserted bounds hold at $\xi=0$.
		
		Lemma~\ref{I:lem:uniform-C2-bubble-map} gives
		\[
		\|D_\xi\Theta_{\bar z}(\xi)-D_\xi\Theta_{\bar z}(0)\|_{\mathrm{op}}
		\le C|\xi|.
		\]
		Choosing $r$ small proves the bounds throughout $|\xi|\le r$.
	\end{proof}

	\begin{lemma}\label{I:lem:preliminary-identification}
		There exist $\eps_1>0$ and $\lambda_0>1$ such that, for two separated
		configurations with internally comparable scales satisfying
		\[
		\lambda_t,\lambda_t'\ge\lambda_0\quad(t\in I'),
		\qquad
		\|U_{p,\lambda,0}-U_{p',\lambda',0}\|_{\mathcal H}<\eps_1,
		\]
		one can relabel the second configuration so that
		\[
		\max_t\lambda_t d_g(p_t,p_t')
		+\max_t\left|\log\frac{\lambda_t'}{\lambda_t}\right|
		\le C\|U_{p,\lambda,0}-U_{p',\lambda',0}\|_{\mathcal H}.
		\]
	\end{lemma}
	
	\begin{proof}
		First consider sequences for which the difference tends to zero in
		$\mathcal H$ and both minimum scales tend to infinity.
		Testing the defining Poisson equations and using invertibility of $A$
		shows that their limiting atomic measures agree:
		\[
		\sum_t\delta_{p_t^*}=\sum_t\delta_{p_t^{\prime *}}.
		\]
		Separation matches the limiting centers after one permutation.
		For a matched pair, set
		\[
		a_k:=\frac{\lambda_t'^k}{\lambda_t^k},\qquad
		d_k:=\lambda_t^k y_{p_t^k}(p_t'^k).
		\]
		Tests of the form $\zeta(\lambda_t^k y_{p_t^k})$ have bounded Dirichlet
		norm. The Poisson equations therefore give the same limits for the
		rescaled densities $\chi_p e^{-\varphi_p}e^{W_{i,p,\lambda}}dV_g$.
		The first density converges to $e^{v_i(y)}dy$, whereas the second has
		the bound
		\[
		a_k^2e^{v_i(a_k(y-d_k))}
		\le\frac{Ca_k^2}{(1+a_k|y-d_k|)^{m_i}}
		\le\frac C{|y-d_k|^2}.
		\]
		This bound and $\int_{\mathbb R^2}e^{v_i}=2\pi\sigma_i$ give the
		locally weak limits
		\[
		a_k^2e^{v_i(a_k(y-d_k))}\,dy\ \rightharpoonup_{\mathrm{loc}}
		\begin{cases}
			0,&a_k\to0\ \text{or }|d_k|\to\infty,\\
			2\pi\sigma_i\delta_d,&a_k\to\infty,\ d_k\to d.
		\end{cases}
		\]
		The normalized chart transitions are asymptotically orthogonal, which
		preserves the radial profile. Both displayed limits contradict
		the smooth positive first density. Thus $a_k$ stays bounded away from
		zero and infinity, and $d_k$ is bounded. Any subsequential limit satisfies  
		\[
		e^{v_i(y)}=a^2e^{v_i(a(y-d))}.
		\]
		Since $rv_i'(r)=-\sum_j a_{ij}\int_0^r s e^{v_j(s)}ds<0$,
		the maximum is unique. Thus $d=0$, and evaluation at zero gives $a=1$.
		Thus the normalized parameter differences tend to zero.
		
		By contradiction, sufficiently small $\eps_1$ and large $\lambda_0$
		place every pair in one fixed small normalized parameter neighborhood.
		With the first parameters as $\bar z$, Taylor's formula and
		Lemmas~\ref{I:lem:uniform-C2-bubble-map} and \ref{I:lem:tangent-gram} yield
		\[
		\|\Theta_{\bar z}(\xi)-\Theta_{\bar z}(0)\|_{\mathcal H}
		\ge c|\xi|-C|\xi|^2\ge\frac c2|\xi|.
		\]
		Since $|\xi|$ is equivalent to the displayed normalized parameter
		distance, the quantitative estimate follows.
	\end{proof}

	\medskip
	\noindent\textbf{Acknowledgements.}~
	During the preparation of this manuscript, the authors used ChatGPT (OpenAI, GPT-6.1 Sol)  to assist with organization, LaTeX consistency, proof
	checking, and the formulation of auxiliary arguments. The authors subsequently reviewed, revised, and verified all content and take full responsibility for the accuracy, originality, and integrity of the manuscript.
	
	\section*{Statements and Declarations}
	\noindent\textbf{Funding.}~
	M. Ahmedou is supported by the Deutsche Forschungsgemeinschaft (DFG) under grant AH~156/2-1. Z. Hu is partially supported by the
	National Natural Science Foundation of China
	(Grant No.~12601383) and the China Postdoctoral Science Foundation
	(Grant No.~2026M793392).  L. Zhang is partially supported by Simons Foundation Grant SFI-MPS-TSM-00013752.

	\medskip
	
	\noindent\textbf{Data availability.}~	 Data sharing is not applicable to this article as no datasets were generated or analysed during the current study.

	\medskip
	
	\noindent\textbf{Competing interests.} ~ The authors declare that there is no conflict of interest.

\end{document}